\documentclass[leqno,final]{siamltex}
\usepackage{graphicx} 
\usepackage{amsmath,amstext,amssymb,bm}
\usepackage{leftidx}
\usepackage{listings}

\usepackage{xcolor} 
\usepackage{soul} 
\usepackage{tikz}
\usetikzlibrary{shapes,arrows}
\usepackage{mathrsfs}
\usepackage{epstopdf}
\usepackage{color}
\usepackage{multirow}
\usepackage{tabularx}
\usepackage[shortlabels]{enumitem}

\numberwithin{equation}{section}

\allowdisplaybreaks[4]

\renewcommand{\div}{\mbox{\rm div\,}}
\newcommand{\curl}{\mbox{\rm curl\,}}

\newcommand{\eps}{\epsilon}

\newcommand{\p}{\partial}

\definecolor{codegreen}{rgb}{0,0.6,0}
\definecolor{codegray}{rgb}{0.5,0.5,0.5}
\definecolor{codepurple}{rgb}{0.58,0,0.82}
\definecolor{backcolour}{rgb}{0.95,0.95,0.92}

\lstdefinestyle{mystyle}{
	backgroundcolor=\color{backcolour},   
	commentstyle=\color{codegreen},
	keywordstyle=\color{magenta},
	numberstyle=\tiny\color{codegray},
	stringstyle=\color{codepurple},
	basicstyle=\ttfamily\footnotesize,
	breakatwhitespace=false,         
	breaklines=true,                 
	captionpos=b,                    
	keepspaces=true,                 
	numbers=left,                    
	numbersep=5pt,                  
	showspaces=false,                
	showstringspaces=false,
	showtabs=false,                  
	tabsize=2
}

\newcommand{\Integer}{\mbox{${\bf Z}$}}
\newcommand{\Real}{\mbox{${\bf R}$}}

\newcommand{\Grad}[1]{\nabla #1}

\newcommand{\Div}[1]{\div\left[#1\right]}

\newcommand{\Curl}[1]{\text{curl}\left[#1\right]}
\newcommand{\Laplacian}[1]{\Delta #1}

\newcommand{\Norm}[2]{\left\|#1\right\|_{#2}}

\newcommand{\SupNorm}[1]{\left|#1\right|_{L^{\infty}}}
\newcommand{\HolderNorm}[2]{\left|#1\right|_{C^{#2}}}
\newcommand{\SobNorm}[2]{\left\|#1\right\|_{H^{#2}}}

\newcommand{\Abs}[1]{\left|#1\right|}

\newcommand{\Angle}[1]{\langle #1 \rangle}

\newcommand{\cL}{\mathcal{L}}

\newcommand{\cP}{\mathcal{P}}

\newcommand{\px}{\partial_x}
\newcommand{\py}{\partial_y}
\newcommand{\pz}{\partial_z}

\newcommand{\sump}{\sum_{p=-\infty}^{\infty}}
\newcommand{\sumq}{\sum_{q=-\infty}^{\infty}}

\newcommand{\be}{\begin{equation}}
	\newcommand{\ee}{\end{equation}}
\newcommand{\bes}{\begin{equation*}}
	\newcommand{\ees}{\end{equation*}}
\newcommand{\bse}{\begin{subequations}}
	\newcommand{\ese}{\end{subequations}}

\newcommand{\void}[1]{}

\newcommand{\epsu}{\epsilon^{(u)}}
\newcommand{\epsv}{\epsilon^{(v)}}
\newcommand{\epsw}{\epsilon^{(w)}}
\newcommand{\epsm}{\epsilon^{(m)}}

\newcommand{\gammau}{\gamma^{(u)}}

\newcommand{\gammaw}{\gamma^{(w)}}
\newcommand{\gammae}{\gamma^{(\bar{\epsilon})}}
\newcommand{\gammam}{\gamma^{(m)}}
\newcommand{\bgamma}{\bar{\gamma}}

\newcommand{\ku}{k^{(u)}}

\newcommand{\km}{k^{(m)}}

\newcommand{\bepsilon}{\bar{\epsilon}}

\newcommand{\tepsilon}{\mathcal{E}}

\begin{document}
	
	\title{A High--Order Perturbation of Envelopes (HOPE) Method for
		Vector Electromagnetic Scattering by Periodic Inhomogeneous Media
		\thanks{D.P.N. gratefully acknowledges support from the National Science Foundation
			through Grant No.~DMS--2111283.}}
	\markboth{DAVID P.\ NICHOLLS AND LIET VO}{HIGH--ORDER PERTURBATION OF ENVELOPES METHOD}
	\author{
		David P.\ Nicholls 
		\thanks{Department of Mathematics, Statistics and Computer Science,
			University of Illinois at Chicago, Chicago, IL 60607, U.S.A. (davidn@uic.edu).}
		\and 
		Liet Vo
		\thanks{Department of Mathematics, Statistics and Computer Science,
			University of Illinois at Chicago, Chicago, IL 60607, U.S.A. (lietvo@uic.edu).}
	}
	
	\maketitle
	
	\begin{abstract}
		The scattering of electromagnetic waves by three--dimensional
		periodic structures is important for many problems of crucial
		scientific and engineering interest. Due to the complexity
		and three--dimensional nature of these waves, the fast, accurate,
		and reliable
		numerical simulations of these are indispensable for engineers and
		scientists alike. For this, High--Order Spectral
		methods are frequently employed and here we describe an algorithm
		in this class. Our approach is perturbative in nature where we
		view the deviation of the permittivity from a constant value as
		the deformation and we pursue regular perturbation theory.
		This work extends our previous contribution regarding
		the Helmholtz equation to the full vector Maxwell equations,
		by providing a rigorous analyticity theory, both in deformation
		size and spatial variable (provided that the permittivity is,
		itself, analytic).
	\end{abstract}
	
	
	
	%
	%
	
	\section{Introduction}
	\label{Sec:Intro}
	
	The scattering of electromagnetic waves by three--dimensional
	periodic structures is important for many problems of crucial
	scientific and engineering interest. Examples abound in areas
	as disparate as surface enhanced spectroscopy \cite{Moskovits85},
	extraordinary optical transmission \cite{ELGTW98},
	cancer therapy \cite{ESHES06}, and surface plasmon resonance (SPR)
	biosensing
	\cite{Homola08,JJJLWO13,LJJOO12,NichollsReitichJohnsonOh14}.
	
	Due to their central role in these nanotechnologies, simulations
	of these waves
	have been conducted with all of the classical numerical algorithms
	for approximating solutions to the relevant
	governing partial differential equations. This includes
	the Finite Difference \cite{Strikwerda04,LeVeque07}, Finite Element
	\cite{Johnson87,Ihlenburg98}, Discontinuous Galerkin
	\cite{HesthavenWarburton08}, Spectral Element \cite{DevilleFischerMund02},
	and Spectral \cite{GottliebOrszag77,ShenTang06,ShenTangWang11} methods.
	While these are tempting choices, due to their \textit{volumetric}
	nature they require a large number of unknowns
	($N = N_x N_y N_z$ for a three--dimensional simulation)
	and mandate the inversion of large, non--symmetric positive definite
	matrices (of dimension $N \times N$). Such properties are still
	an object of current study (for instance, see
	\cite{ErnstGander12,MoiolaSpence14}).
	
	For the specific application of SPR sensors which we have in mind
	in the current contribution, their pervasiveness stems from two 
	properties of an SPR, namely its extremely strong and sensitive
	response. More specifically, over
	the range of tens of nanometers in incident wavelength, the reflected
	energy can fall from nearly 100~\% by a factor of 10 or even 100
	before returning to almost 100~\%. Obviously, to approximate such
	a structure with the required accuracy, the numerical algorithm should
	produce high fidelity results in a rapid and robust manner.
	For this reason we will focus upon High--Order Spectral (HOS) methods
	\cite{GottliebOrszag77,ShenTang06,ShenTangWang11} which can deliver
	precisely this behavior.
	
	Returning to the classical approaches listed above, for the problem of
	scattering by homogeneous layers (which is one important avenue to
	generating SPRs) it is clearly unnecessary
	to discretize the bulk of each layer and state--of--the--art solvers
	seek interfacial unknowns with the knowledge that information
	inside the layers can readily be computed from appropriate integral
	formulas. Boundary element (BEM) \cite{SauterSchwab11} and boundary
	integral (BIM) \cite{ColtonKress13,Kress14} methods are two such approaches
	and can produce spectrally accurate solutions in a
	fraction of the time of their volumetric competitors.
	
	In previous work \cite{Nicholls19b} the authors investigated a novel
	algorithm very much in the spirit of these HOS approaches, and inspired
	by the ``High--Order Perturbation of Surfaces'' (HOPS) algorithms which
	have proven to be so appealing for layered media. A HOPS scheme is one
	which views the layer interfaces as perturbations of flat ones and then
	makes recursive corrections to the scattering returns from this classical,
	exactly solvable, configuration \cite{Yeh05}. By contrast, our
	new ``High--Order Perturbation of Envelopes'' (HOPE) scheme considers
	more a general permittivity function, $\epsilon(x,y,z)$, which does not
	necessarily have layered structure. We followed the lead of
	Feng, Lin, and Lorton \cite{FengLinLorton15,FengLinLorton16} and 
	adopted a perturbative philosophy (much like a HOPS algorithm) by
	viewing the permittivity as a perturbation of a trivial one, e.g.,
	\bes
	\epsilon(x,y,z) = \bepsilon - \delta \tepsilon(x,y,z),
	\quad
	\bepsilon \in \Real,
	\quad
	\tepsilon(x+d_x,y+d_y,z) = \tepsilon(x,y,z),
	\ees
	where $\tepsilon$ is a permittivity ``envelope.''
	In this previous contribution we focused upon the two--dimensional
	scalar problems governing electromagnetic radiation in
	Transverse Electric (TE) and Transverse Magnetic (TM) polarizations.
	This new approach has computational advantages over
	volumetric solvers in certain configurations (e.g., where the support
	of $\tepsilon$ is small or where the set on which $\tepsilon$
	significantly changes is small). A particular choice we pursued
	was an approximate indicator function which
	is nearly zero/unity to denote the absence/presence
	of a material.
	
	Among the contributions of \cite{Nicholls19b} was a new, extensive, 
	and rigorous analysis. More specifically,
	we proved not only that the domain of analyticity of the scattered
	field in $\delta$ can be extended to a neighborhood of the
	\textit{entire} real axis (up to topological obstruction),
	but also that this field is \textit{jointly} analytic
	in parametric and spatial variables provided that
	$\tepsilon(x,y,z)$ is spatially analytic.
	In the current paper we extend some of these results to the
	three--dimensional vector electromagnetic case governed by
	the full time--harmonic Maxwell equations. More specifically,
	we show that the scattered field is analytic as a function
	of $\delta$ and \textit{jointly} analytic in both parametric
	and spatial variables if $\tepsilon(x,y,z)$ is spatially analytic.
	We delay for future consideration the issue of the analytic 
	continuation of our
	results to perturbations $\delta$ of arbitrary (real) size. This
	requires an analysis of the variable coefficient Maxwell equations
	\cite{BaoLi22} which is more subtle than what we present here, 
	however, 
	we believe that this is a result which can be established with our
	current framework.
	
	The rest of the paper is organized as follows. In \S~\ref{Sec:Govern}
	we recall the governing equations complete with a discussion of
	transparent boundary conditions in \S~\ref{Sec:TransBC}. We describe the
	HOPE algorithm in \S~\ref{Sec:HOPE} and begin our theoretical
	developments with a description of the relevant function spaces
	in \S~\ref{Sec:Func}. We state and prove our results on parametric
	analyticity in \S~\ref{Sec:Anal} and joint parametric/spatial analyticity
	in \S~\ref{Sec:JointAnal}.
	In Appendix~\ref{Sec:EllEst:Proof} we prove our new elliptic
	estimate, and in Appendix~\ref{Sec:JointAnal:Proof} we establish the
	joint analyticity result in the special case of lateral smoothness.
	
	%
	%
	
	\section{Governing Equations}
	\label{Sec:Govern}
	
	We consider materials modeled by the time--harmonic Maxwell 
	equations in three dimensions with a constant permeability
	$\mu = \mu_0$ and no currents or sources,
	\be
	\label{Eqn:Maxwell:TimeHarmonic}
	\Curl{E} - i \omega \mu_0 H = 0,
	\quad
	\Curl{H} + i \omega \eps E = 0,
	\ee
	where $(E, H)$ are the electric and magnetic vector fields, and we have
	factored out time dependence of the form $\exp(-i \omega t)$ \cite{BaoLi22}.
	The permittivity $\epsilon(x,y,z)$ is biperiodic with periods $d_x$ and
	$d_y$, and is specified by
	\bes
	\eps(x,y,x) = \begin{cases}
		\epsu, & z > h, \\
		\epsv(x,y,z), & -h < z < h, \\
		\epsw, & z < -h,
	\end{cases}
	\ees
	where $\epsu, \epsw \in \Real^+$, and $\epsv(x+d_x,y+d_y,z) = \epsv(x,y,z)$,
	and
	\bes
	\lim_{z \rightarrow h-} \epsv(x,y,z) = \epsu,
	\quad
	\lim_{z \rightarrow (-h)+} \epsv(x,y,z) = \epsw.
	\ees
	Using the permittivity of vacuum, $\epsilon_0$, we can define
	\bes
	k^2_0 = \omega^2 \epsilon_0 \mu_0 = \frac{\omega^2}{c_0^2},
	\quad
	(\km)^2 = \epsm k^2_0,
	\quad
	m \in \{u, w\},
	\ees
	and $c_0 = 1/\sqrt{\epsilon_0 \mu_0}$ is the speed of light in the vacuum.
	
	This structure is illuminated from above by plane--wave incident 
	radiation of the form 
	\begin{align*}
		E^{\text{inc}}(x,y,z) & = A \exp(i \alpha x + i \beta y - i \gammau z), \\
		H^{\text{inc}}(x,y,z) & = B \exp(i \alpha x + i \beta y - i \gammau z), 
	\end{align*}
	where
	\bes
	A \cdot \kappa = 0,
	\quad
	B = \frac{1}{\omega \mu_0} \kappa \times A,
	\quad
	\Abs{A} = \Abs{B} = 1,
	\ees
	and
	\bes
	\kappa = \begin{pmatrix} \alpha \\ \beta \\ -\gammau \end{pmatrix}
	= \ku \begin{pmatrix} \sin(\theta) \cos(\phi) \\
		\sin(\theta) \sin(\phi) \\ -\cos(\theta) \end{pmatrix},
	\ees
	where $(\theta,\phi)$ are the angles of incidence.
	
	%
	%
	
	\section{Transparent Boundary Conditions}
	\label{Sec:TransBC}
	
	Following our previous work \cite{Nicholls19b} we aim to
	both rigorously specify the appropriate far--field boundary
	conditions and reduce the infinite domain to one of finite
	size. Conveniently, a variant of the scheme presented in
	Bao \& Li \cite{BaoLi22} accomplishes both. In the upper
	domain $\{ z > h \}$ we seek a solution as the sum of
	the incident radiation and an upward propagating (reflected) component,
	e.g.,
	\begin{align*}
		E & = E^{\text{inc}} + E^{\text{refl}} \\
		& = A \exp(i \alpha x + i \beta y - i \gammau z) 
		+ c_{0,0} \exp(i \alpha x + i \beta y + i \gammau (z-h)) \\
		& \quad 
		+ \sum_{(p,q) \neq (0,0)} \hat{u}_{p,q} \exp(i \alpha_p x 
		+ i \beta_q y + i \gammau_{p,q} (z-h)),
	\end{align*}
	\cite{Petit80,Yeh05} where
	\begin{gather*}
		\alpha_p = \alpha + (2 \pi/d_x) p,
		\quad
		\beta_q = \beta + (2 \pi/d_y) q,
		\\
		{\gammam_{p,q} := \begin{cases}
				\sqrt{\epsilon^{(m)} k_0^2 - \alpha_p^2 - \beta_q^2},
				& \alpha_p^2 + \beta_q^2 \leq \epsilon^{(m)} k_0^2, \\
				i \sqrt{\alpha_p^2 + \beta_q^2 - \epsilon^{(m)} k_0^2},
				& \alpha_p^2 + \beta_q^2 > \epsilon^{(m)} k_0^2,
		\end{cases}}
		\quad
		m \in \{ u, w \}.
	\end{gather*}
	If we set $c_{0,0} = \hat{u}_{0,0} - A \exp(-i \gammau h)$, implying
	$\hat{u}_{0,0} = c_{0,0} + A \exp(-i \gammau h)$, then
	\begin{align*}
		E & = A \exp(i \alpha x + i \beta y - i \gammau z) 
		- A \exp(i \alpha x + i \beta y + i \gammau (z-2h)) \\
		& \quad
		+ \sump \sumq \hat{u}_{p,q} 
		\exp(i \alpha_p x + i \beta_q y + i \gammau_{p,q} (z-h)),
	\end{align*}
	and $E(x,y,h) = u(x,y)$. It is a simple matter to show that
	\begin{align*}
		\pz E 
		& = (-i \gammau) A \exp(i \alpha x + i \beta y - i \gammau z) \\
		& \quad 
		- (i \gammau) A \exp(i \alpha x + i \beta y + i \gammau (z-2h)) \\
		& \quad
		+ \sump \sumq (i \gammau_{p,q}) \hat{u}_{p,q} 
		\exp(i \alpha_p x + i \beta_q y + i \gammau_{p,q} (z-h)),
	\end{align*}
	so that
	\begin{align*}
		-\pz E(x,y,h) 
		& = (i \gammau) A \exp(i \alpha x + i \beta y - i \gammau h) \\
		& \quad
		+ (i \gammau) A \exp(i \alpha x + i \beta y + i \gammau (-h)) \\
		& \quad
		+ \sump \sumq (-i \gammau_{p,q}) \hat{u}_{p,q} 
		\exp(i \alpha_p x + i \beta_q y).
	\end{align*}
	If we define the function
	\be
	\label{Eqn:phi:Def}
	\phi(x,y) := \left( 2 i \gammau \exp(-i \gammau h) \right) A 
	\exp(i \alpha x + i \beta y),
	\ee
	and the order--one Fourier multiplier (the externally directed
	Dirichlet--Neumann operator for the Maxwell equation on $\{ z > h \}$)
	\bes
	T_u[ \psi ] := \sump \sumq (-i \gammau_{p,q}) 
	\hat{\psi}_{p,q} \exp(i \alpha_p x + i \beta_q y),
	\ees
	then we see that we can express the Upward Propagating Condition (UPC)
	\cite{ArensHab} exactly with the boundary condition
	\bes
	-\pz E(x,y,h) - T_u[ E(x,y,h) ] = \phi(x,y).
	\ees
	
	In a similar fashion, in $\{ z < -h \}$ we seek a solution
	which is purely downward propagating (transmitted)
	\bes
	E = E^{\text{trans}}
	= \sump \sumq \hat{w}_{p,q} 
	\exp(i \alpha_p x + i \beta_q y - i \gammaw_{p,q} (z+h)),
	\ees
	\cite{Petit80,Yeh05}.
	Clearly $E(x,y,-h) = w(x,y)$ and, with the calculation
	\bes
	\pz E(x,y,-h) = \sump \sumq (-i \gammaw_{p,q}) \hat{w}_{p,q}
	\exp(i \alpha_p x + i \beta_q y),
	\ees
	and the analogous order--one Fourier multiplier (again, the
	externally directed Dirichlet--Neumann operator for the Helmholtz 
	equation on $\{ z < -h \}$)
	\bes
	T_w[ \psi ] := \sump \sumq (-i \gammaw_{p,q}) \hat{\psi}_{p,q}
	\exp(i \alpha_p x + i \beta_q y),
	\ees
	we can state the Downward Propagating Condition (DPC) \cite{ArensHab}
	transparently using
	\bes
	\pz E(x,y,-h) - T_w[ E(x,y,-h) ] = 0.
	\ees
	
	Eliminating the magnetic field from \eqref{Eqn:Maxwell:TimeHarmonic}
	and gathering our full set of governing equations we find
	the following problem to solve.
	\bse
	\label{Eqn:Max}
	\begin{align}
		& \Curl{ \Curl{ E } } - \epsv k_0^2 E = 0,
		&& \text{in $S_v$}, \label{Eqn:Max:a} \\
		& -\pz E - T_u[ E ] = \phi, && \text{at $\Gamma_h$},
		\label{Eqn:Max:b} \\
		& \pz E - T_w[ E ] = 0, && \text{at $\Gamma_{-h}$},
		\label{Eqn:Max:c} \\
		& E(x+d_x,y+d_y,z) = \exp(i \alpha d_x + i \beta d_y) E(x,y,z), 
		\label{Eqn:Max:d}
	\end{align}
	\ese
	where
	\bes
	S_v := (0,d_x) \times (0,d_y) \times (-h,h),
	\quad
	\Gamma_{\pm h} := (0,d_x) \times (0,d_y) \times \{ z = \pm h \}.
	\ees
	
	
	%
	%
	
	\section{A High--Order Perturbation of Envelopes Method}
	\label{Sec:HOPE}
	
	Following the lead of our previous work \cite{Nicholls19b} we
	do not pursue the solution of \eqref{Eqn:Max} by a classical
	volumetric approach, but rather a perturbative one where 
	we think of our configuration as a small deviation from a trivial
	structure,
	\bes
	\epsv(x,y,z) = \bepsilon \left( 1 - \delta \tepsilon(x,y,z) \right)
	= \bepsilon - \delta (\bepsilon \tepsilon(x,y,z)),
	\ees
	where $\bepsilon \in \Real^+$ is a constant, and $\delta \ll 1$.
	In our previous work on the Helmholtz equation in either TE or TM
	polarization, we showed that if $\tepsilon(x,z)$ is smooth enough
	then the transverse components of $E$ or $H$ depend analytically
	upon $\delta$. In light of this we posit that, for the full Maxwell
	equations which we consider here, the field $E = E(x,y,z;\delta)$
	depends analytically upon $\delta$ so that
	\be
	\label{Eqn:EExp}
	E = E(x,y,z;\delta) = \sum_{\ell=0}^{\infty} E_{\ell}(x,y,z) \delta^{\ell},
	\ee
	converges strongly in a function space. It is not difficult to see that
	these $E_{\ell}$ must satisfy
	\bse
	\label{Eqn:Max:ell}
	\begin{align}
		& \Curl{ \Curl{ E_{\ell} } } - \bepsilon k_0^2 E_{\ell} 
		= \bepsilon k_0^2 \tepsilon E_{\ell-1},
		&& \text{in $S_v$}, \label{Eqn:Max:ell:a} \\
		& -\pz E_{\ell} - T_u[ E_{\ell} ] = \delta_{\ell,0} \phi,
		&& \text{at $\Gamma_h$}, \label{Eqn:Max:ell:b} \\
		& \pz E_{\ell} - T_w[ E_{\ell} ] = 0, 
		&& \text{at $\Gamma_{-h}$}, \label{Eqn:Max:ell:c} \\
		& E_{\ell}(x+d_x,y+d_y,z) 
		= \exp(i \alpha d_x + i \beta d_y) E_{\ell}(x,y,z),
		\label{Eqn:Max:ell:d}
	\end{align}
	\ese
	where $\delta_{\ell,0}$ is the Kronecker delta function. It is
	easy to see that
	\be
	\label{Eqn:E0}
	E_0(x,y,z) = A e^{i \alpha x + i \beta y - i \bgamma z},
	\quad
	\alpha^2 + \beta^2 + (\bgamma)^2 = \bepsilon k_0^2,
	\ee
	and our HOPE scheme can be viewed as computing corrections to
	this by
	\bes
	E(x,y,z) = A e^{i \alpha x + i \beta y - i \bgamma z}
	+ \sum_{\ell=1}^{\infty} E_{\ell}(x,y,z) \delta^{\ell}.
	\ees
	We note that a natural numerical method would seek approximations
	of the $E_{\ell}$ for $1 \leq \ell \leq L$ and simulate
	\bes
	E(x,y,z;\delta) \approx E^L(x,y,z;\delta)
	:= \sum_{\ell=0}^{L} E_{\ell}(x,y,z) \delta^{\ell}.
	\ees
	
\begin{figure}[hbt]
		\begin{center}
			\includegraphics[width=0.45\textwidth]{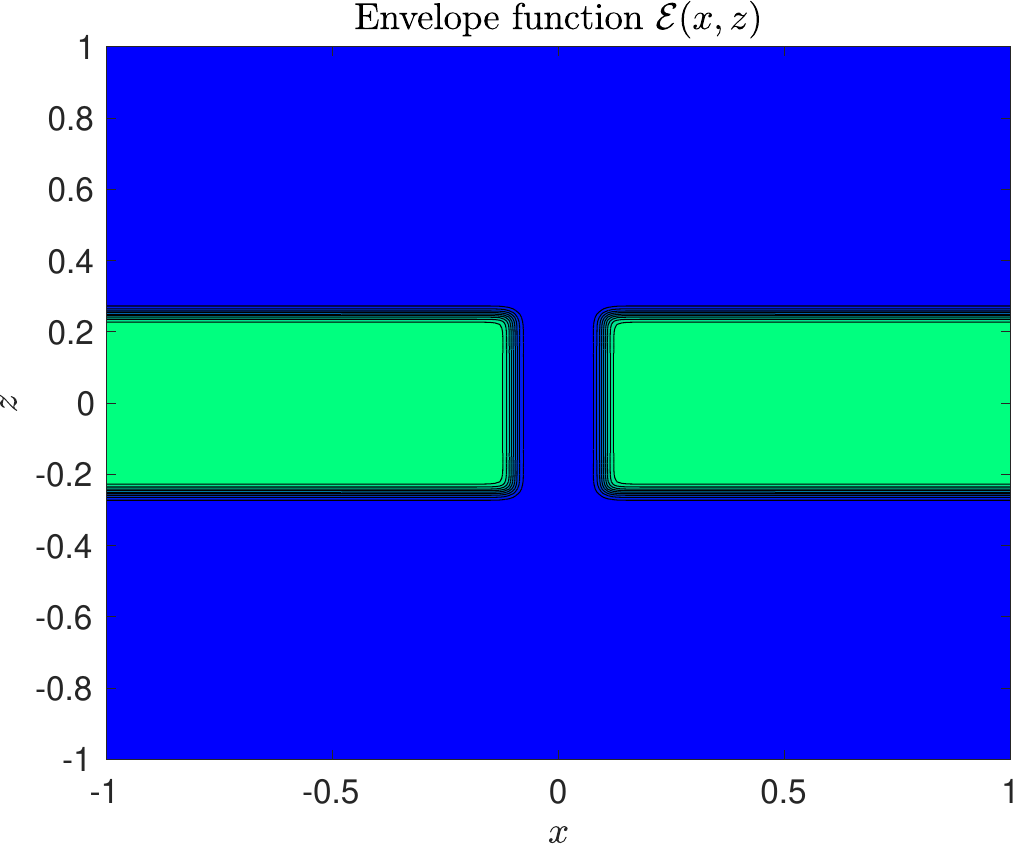}
			\quad
			\includegraphics[width=0.45\textwidth]{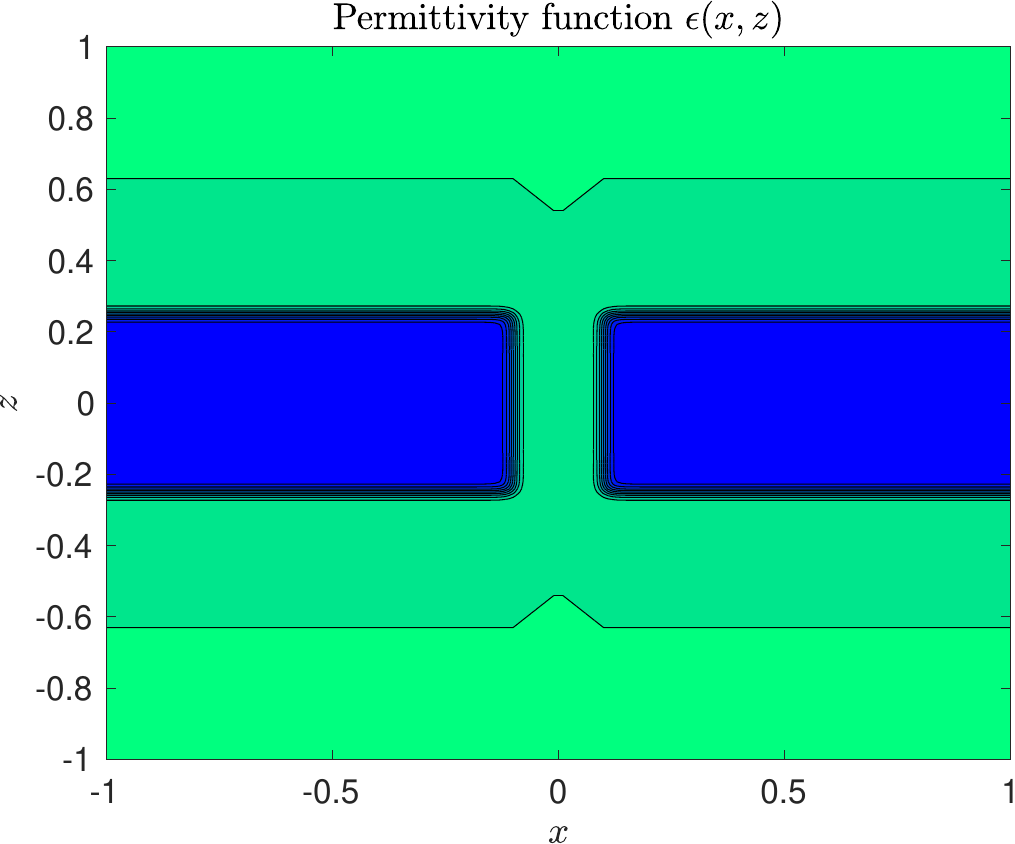}
			\caption{Contour plots of $\tepsilon(x,z)$ (left) and 
				$\epsv(x,z)$ (right).}
			\label{Fig:E:Perm}
		\end{center}
	\end{figure}

	There are many possibilities for the envelope function
	$\tepsilon(x,y,z)$ and each leads to a slightly different
	perturbation approach. For instance, consider the function
	\bes
	\Phi_{a,b}(z) := \frac{\tanh(w (z-a)) - \tanh(w (z-b))}{2},
	\ees
	with sharpness parameter $w$, which is effectively zero outside
	the interval $(a,b)$ while being essentially one inside
	$(a,b)$, c.f. \cite{Nicholls19b}. We can 
	approximate a slab of material with thickness $2d$ and a gap of
	width $2g$ in vacuum by selecting
	\bes
	\bepsilon = 1,
	\quad
	\tepsilon(x,y,z) = \Phi_{-d,d}(z) 
	\left\{ 1 - \Phi_{-g,g}(x) \right\},
	\quad
	\delta = 1.
	\ees
	See Figure~\ref{Fig:E:Perm} with the choices $d=1/4$, $g=1/10$,
	and $w=50$ on the cell $[-1,1] \times [-1,1]$.
	%
	%
	
	%
	%
	
	\section{Function spaces}
	\label{Sec:Func}
	
	In this section, we present function spaces and theoretical
	notions that are necessary for our analysis later. We point
	out that due to the vectorial nature of the Maxwell equations,
	these are extensions of the results we utilized in our previous
	study \cite{Nicholls19b}. For any real number $s \geq 0$,
	we have the classical interfacial quasiperiodic $L^2$ Sobolev norm
	\begin{gather*}
		\SobNorm{v}{s}^2 := \sump \sumq \Angle{(p,q)}^{2 s}
		\Abs{\hat{v}_{p,q}}^2, \\
		\Abs{\hat{v}_{p,q}}^2 = \Abs{\hat{v}^x_{p,q}}^2
		+ \Abs{\hat{v}^y_{p,q}}^2 + \Abs{\hat{v}^z_{p,q}}^2,
	\end{gather*}
	where
	\bes
	\Angle{(p,q)}^2 := 1 + \Abs{p}^2 + \Abs{q}^2,
	\quad
	\hat{v}^m_{p,q} := \frac{1}{d_x d_y} \int_0^{d_x} 
	\int_0^{d_y} v^m(x,y) e^{-i \alpha_p x - i \beta_q y}
	\; dx \; dy,
	\ees
	for $m \in \{x, y, z\}$.
	From this we define the interfacial quasiperiodic Sobolev space
	\cite{Kress14}
	\bes
	H^s(\Gamma) = \left\{ v(x,y) \in (L^2(\Gamma))^3\ |\ 
	\SobNorm{v}{s} < \infty \right\},
	\ees
	where
	\bes
	\Gamma := (0,d_x) \times (0,d_y).
	\ees
	In addition, we mention that the dual space of $H^s$, $H^{-s}$, can be
	defined by the norm above with a negative index. We also recall
	the space of $s$-times continuously differentiable functions with
	H\"older norm
	\bes
	\HolderNorm{v}{s} := \max_{0 \leq \ell + r \leq s} \max_{m \in \{x, y, z\}}
	\SupNorm{\partial_x^{\ell} \partial_y^r v^m}.
	\ees
	Finally, we define the volumetric laterally quasiperiodic Sobolev 
	space as
	\bes
	H^s(S_v) = \left\{ u(x,y,z) \in (L^2(S_v))^3\ |\ 
	\SobNorm{u}{s} < \infty \right\},
	\ees
	where 
	\bes
	\SobNorm{u}{s}^2
	:= \sum_{j=0}^{s} \sump \sumq \Angle{(p,q)}^{2 (s-j)} 
	\int_{-h}^h \Abs{\partial_z^j \hat{u}_{p,q}(z)}^2 \; dz.
	\ees
	We will also have need of the ``H--div'' space
	\bes
	H^s(\div,S_v) = \left\{ u(x,y,z) \in H^s(S_v)\ |\
	\Div{u} \in H^s(S_v) \right\}.
	\ees
	
	We close with an essential result \cite{Evans10,NichollsReitich99} required
	for our later proofs.
	\begin{lemma}
		Let $s \geq 0$ be an integer and $g \in C^s(D)$, $w \in H^s(D)$, where $D$ is a 
		subset of $\Real^m$, then $g w \in H^s(D)$ and
		\bes
		\SobNorm{g w}{s} \leq \tilde{M}(m,s,D) \HolderNorm{g}{s} \SobNorm{w}{s},
		\ees
		where $\tilde{M}$ is some positive universal constant.
	\end{lemma}
	
	Furthermore, we recall the following elementary result
	\cite{NichollsReitich00b,Nicholls19b}.
	\begin{lemma}
		\label{Lemma:S}
		Let $s\geq0$ be an integer, then there exists a constant $S>0$ such that
		\bes
		\sum_{j=0}^s \frac{(s+1)^2}{(s-j+1)^2(j+1)^2} < S,
		\quad
		\sum_{j=0}^s \sum_{r=0}^j \frac{(s+1)^2}{(s-j+1)^2 (j-r+1)^2 (r+1)^2} < S^2.
		\ees
	\end{lemma}
	
	%
	%
	
	\section{Analyticity}
	\label{Sec:Anal}
	
	At this point we are in a position to extend our previous results
	\cite{Nicholls19b} by demonstrating the analytic dependence
	of the full electric field, $E = E(x,y,z;\delta)$, upon $\delta$, 
	sufficiently small. More specifically, we show that the expansion 
	\eqref{Eqn:EExp} converges strongly in an appropriate function space.
	
	For this we require an elliptic estimate for our inductive proof 
	which is established in Appendix~\ref{Sec:EllEst:Proof}. For future 
	convenience, we define the following differential operator associated
	to the Maxwell system
	\bes
	\cL_0 := \curl \curl - k^2_0 \bepsilon.
	\ees
	As is well known \cite{BaoLi22}, the issue of \textit{uniqueness} of 
	solutions to the Maxwell problem
	\bse
	\label{Eqn:Max:Uniqueness}
	\begin{align}
		& \cL_0 V = 0, && \text{in $S_v$}, \\
		& -\pz V - T_u[V] = 0, && \text{at $\Gamma_h$}, \\
		& \pz V - T_w[V] = 0, && \text{at $\Gamma_{-h}$}, \\
		& V(x+d_x,y+d_y,z) = \exp(i \alpha d_x + i \beta d_y) V(x,y,z),
	\end{align}
	\ese
	c.f.\ \eqref{Eqn:Max:ell}, which should have only the \textit{trivial} 
	solution $V \equiv 0$, is a subtle one and certain illuminating
	frequencies $\omega$ will induce non--uniqueness in some configurations.
	Unfortunately a precise characterization of the set of forbidden
	frequencies is elusive and all that is known is that it is countable
	and accumulates at infinity \cite{BaoLi22}. To accommodate this
	state of affairs we define the set of permissible configurations
	\be
	\label{Eqn:cP}
	\cP := \left\{ (\omega,\bepsilon)\ |\ 
	V \equiv 0\ \text{is the unique solution of
		\eqref{Eqn:Max:Uniqueness}} \right\}.
	\ee
	
	With this we can now state the following fundamental elliptic
	regularity result.
	\begin{theorem}
		\label{Thm:EllEst}
		Given any integer $s \geq0$, if
		$(\omega,\bepsilon) \in \cP$,
		$F \in H^{s}(S_v) \cap H^{s+1}(\div,S_v)$,
		$Q \in H^{s+1/2}(\Gamma)$, and
		$R \in H^{s+1/2}(\Gamma)$, 
		then there exists a unique solution of
		\bse
		\label{Eqn:EllEst}
		\begin{align}
			& \cL_0 v = F,  && \text{in $S_v$}, \\
			& -\pz v - T_u[v] = Q && \text{at $\Gamma_h$}, \\
			& \pz v - T_w[v] = R && \text{at $\Gamma_{-h}$}, \\
			& v(x+d_x,y+d_y,z) = e^{i \alpha d_x + i \beta d_y} v(x,y,z),
		\end{align}
		\ese
		satisfying
		\be
		\label{Eqn:EllEst:Est}
		\SobNorm{v}{s+2} \leq C_e \left( \SobNorm{F}{s} 
		+ \SobNorm{\Div{F}}{s+1}
		+ \SobNorm{Q}{s+1/2} + \SobNorm{R}{s+1/2} \right),
		\ee
		where $C_e > 0$ is a universal constant.
	\end{theorem}
	
	We can now prove our analyticity result.
	\begin{theorem}
		\label{Thm:Anal}
		Given any integer $s \geq0$, if
		$(\omega,\bepsilon) \in \cP$,
		and $\tepsilon \in C^{s+2}(S_v)$ then
		the series \eqref{Eqn:EExp} converges strongly.
		More precisely,
		\be
		\label{Eqn:Anal}
		\SobNorm{E_{\ell}}{s+2} \leq K B^{\ell},
		\quad
		\forall\ \ell \geq 0,
		\ee
		for some universal constants $K, B >0$.
	\end{theorem}
	\begin{proof}
		We prove the estimate \eqref{Eqn:Anal} by induction. For $\ell = 0$
		the system \eqref{Eqn:Max:ell} can be written as
		\begin{align*}
			& \cL_0 E_0 = 0, && \text{in $S_v$}, \\
			& -\pz E_0 - T_u [E_0] = \phi, && \text{at $\Gamma_h$}, \\
			& \pz E_0 - T_w[E_0] = 0 && \text{at $\Gamma_{-h}$}, \\
			& E_0(x+d_x,y+d_y,z) = e^{i \alpha d_x + i \beta d_y} E_0(x,y,z).
		\end{align*}
		So, we can apply Theorem~\ref{Thm:EllEst} with $F \equiv 0$,
		$Q = \phi$, and $R \equiv 0$ to obtain
		\bes
		\SobNorm{E_0}{s+2} \leq C_e \SobNorm{\phi}{s+1/2} =: K.
		\ees
		Now, we assume that \eqref{Eqn:Anal} is true for all $\ell <L$ and
		apply Theorem~\ref{Thm:EllEst} to the system \eqref{Eqn:Max:ell}
		for $E_L$ with $F_L = \bepsilon k_0^2 \tepsilon E_{L-1}$ and
		$Q_L \equiv R_L \equiv 0$. This gives
		\begin{align*}
			\SobNorm{E_L}{s+2}
			& \leq C_e \left( \SobNorm{\bepsilon k_0^2 \tepsilon E_{L-1}}{s}
			+ \SobNorm{\Div{\bepsilon k_0^2 \tepsilon E_{L-1}}}{s+1} \right) \\
			& \leq C_e k_0^2 \Abs{\bepsilon} \SobNorm{ \tepsilon E_{L-1}}{s+2} \\
			& \leq C_e k_0^2 \Abs{\bepsilon} \tilde{M} \HolderNorm{\tepsilon}{s+2}
			\SobNorm{E_{L-1}}{s+2} \\
			& \leq C_e k_0^2 \Abs{\bepsilon} \tilde{M} \HolderNorm{\tepsilon}{s+2}
			K B^{L-1} \\
			& \leq K B^L,
		\end{align*}
		provided that
		\bes
		B > C_e k_0^2 \Abs{\bepsilon} \tilde{M} \HolderNorm{\tepsilon}{s+2}.
		\ees
	\end{proof}
	
	From this we can derive the exponential order of convergence
	of the HOPE method. More precisely, defining the $L$--th partial sum
	of \eqref{Eqn:EExp}
	\bes
	E^L(x,y,z) := \sum_{\ell=0}^{L} E_{\ell}(x,y,z) \delta^{\ell},
	\ees
	we obtain the following error estimate for the HOPE method.
	\begin{theorem}
		If $(\omega,\bepsilon) \in \cP$ and $E$ is the unique solution of \eqref{Eqn:Max}, 
		under the assumptions of Theorem~\ref{Thm:Anal} we have the estimate
		\bes
		\SobNorm{E - E^L}{s+2} \lesssim (B \delta)^{L+1},
		\ees
		for some constant
		$B > C_e k_0^2 \Abs{\bepsilon} \tilde{M} \HolderNorm{\tepsilon}{s+2}$
		provided that $\Abs{\delta} < 1/B$.
	\end{theorem}
	\begin{proof}
		Since
		\bes
		E(x,y,z) - E^L(x,y,z) 
		= \sum_{\ell = L+1}^{\infty} E_{\ell}(x,y,z) \delta^{\ell},
		\ees
		we have, by Theorem~\ref{Thm:Anal},
		\bes
		\SobNorm{E - E^L}{s+2}
		\leq \sum_{\ell = L+1}^{\infty} \SobNorm{E_{\ell}}{s+2}
		\leq \sum_{\ell = L+1}^{\infty} K B^{\ell} \delta^{\ell}.
		\ees
		By gathering terms and re--indexing we have
		\bes
		\SobNorm{E - E^L}{s+2}
		\leq K (B \delta)^{L+1} \sum_{\ell=0}^{\infty} (B \delta)^{\ell}
		\leq K \frac{(B \delta)^{L+1}}{1 - (B \delta)},
		\ees
		for $\Abs{\delta B} < 1$, where we have used the elementary fact that
		\bes
		\sum_{\ell=0}^{\infty} \alpha^{\ell} = \frac{1}{1 - \alpha}
		\ees
		provided that $\Abs{\alpha} < 1$.
	\end{proof}
	
	%
	%
	
	\section{Joint Analyticity}
	\label{Sec:JointAnal}
	
	We close our theoretical developments with a result on \textit{joint}
	analyticity of the scattered field with respect to not only the perturbation
	parameter $\delta$, but also the spatial coordinates $(x,y,z)$. As with
	the analogous result for the Helmholtz equation found in \cite{Nicholls19b},
	this will require analyticity of the permittivity envelope function
	$\tepsilon(x,y,z)$. More precisely, we will show that the $E_{\ell}$
	from \eqref{Eqn:EExp} satisfy conditions analogous to those in the following
	definition of analyticity.
	\begin{definition}
		Given an integer $m \geq 0$, if the functions $f = f(x,y)$ and 
		$\tepsilon = \tepsilon(x,y,z)$ are \textit{real analytic}
		and satisfy the following estimates
		\begin{align*}
			\HolderNorm{\frac{\px^r \py^t}{(r+t)!} f}{m} 
			& \leq C_f \frac{\eta^r}{(r+1)^2} \frac{\theta^t}{(t+1)^2}, \\
			\HolderNorm{\frac{\px^r \py^t \pz^s}{(r+t+s)!} \tepsilon}{m} 
			& \leq C_{\tepsilon} \frac{\eta^r}{(r+1)^2} \frac{\theta^t}{(t+1)^2}
			\frac{\zeta^s}{(s+1)^2},
		\end{align*}
		for all $r, t, s \geq 0$,
		for some constants $C_f, C_{\tepsilon}, \eta, \theta, \zeta > 0$,
		then $f \in C^{\omega}_m(\Gamma)$ and
		$\tepsilon \in C^{\omega}_m(S_v)$.
	\end{definition}
	
	The space $C^{\omega}_m$ is the space of real analytic functions
	with radius of convergence (specified by $\eta$, $\theta$, and
	$\zeta$) measured in the $C^m$ norm.
	It is clear that the incident radiation function $\phi$,
	\eqref{Eqn:phi:Def}, is jointly
	analytic in $x$ and $y$ as we now explicitly state.
	\begin{lemma}
		The function
		\bes
		\phi(x,y) := \left( 2 i \gammau \exp(-i \gammau h) \right) A 
		\exp(i \alpha x + i \beta y),
		\ees
		is real analytic and satisfies
		\bes
		\SobNorm{ \frac{\px^r \py^t}{(r+t)!} \phi }{1/2} 
		\leq C_{\phi} \frac{\eta^r}{(r+1)^2} \frac{\theta^t}{(t+1)^2},
		\ees
		for all $r, t \geq 0$, for some constants $C_{\phi}, \eta, \theta > 0$.
	\end{lemma}
	
	Now we present the fundamental elliptic estimate which is required in our forthcoming proof. It is proven in 
	Appendix~\ref{Sec:JointAnal:Proof}.
	\begin{theorem}
		\label{Thm:JointAnal:EllEst}
		Given any integer $m \geq 0$, if
		$(\omega,\bepsilon) \in \cP$,
		$F(x,y,z) \in C^{\omega}_m(S_v)$
		such that
		\bes
		\max \left\{
		\SobNorm{ \frac{\px^r \py^t \pz^s}{(r+t+s)!} F }{0},
		\SobNorm{ \frac{\px^r \py^t \pz^s}{(r+t+s)!} \Div{F} }{0}
		\right\}
		\leq C_F \frac{\eta^r}{(r+1)^2} \frac{\theta^t}{(t+1)^2}
		\frac{\zeta^s}{(s+1)^2},
		\ees
		for all $r, t, s \geq 0$ and for some constants 
		$C_F, \eta, \theta, \zeta > 0$, and
		$Q, R \in C^{\omega}_m(\Gamma))$ satisfying 
		\begin{align*}
			& \SobNorm{ \frac{\px^r \py^t}{(r+t)!} Q }{1/2} 
			\leq C_Q \frac{\eta^r}{(r+1)^2} \frac{\theta^t}{(t+1)^2}, \\
			& \SobNorm{\frac{\px^r \py^t}{(r+t)!} R }{1/2} 
			\leq C_R \frac{\eta^r}{(r+1)^2} \frac{\theta^t}{(t+1)^2},
		\end{align*}
		for all $r, t \geq 0$ and some constants $C_R, C_Q > 0$.
		Then, there exists a unique solution
		$v \in C^{\omega}_m(S_v)$ of 
		\begin{align*}
			& \cL_0 v = F, && \text{in $S_v$}, \\
			& -\pz v - T_u [v] = Q, && \text{at $\Gamma_h$}, \\
			& \pz v - T_w[v] = R, && \text{at $\Gamma_{-h}$}, \\
			& v(x+d_x,y+d_y, z) = e^{i \alpha d_x + i \beta d_y} v(x,y,z),
		\end{align*}
		satisfying 
		\be
		\label{Eqn:JointAnal:EllEst:Estimate}
		\SobNorm{ \frac{\px^r \py^t \pz^s}{(r+t+s)!} v }{2} 
		\leq \underline{C}_e \frac{\eta^r}{(r+1)^2}
		\frac{\theta^t}{(t+1)^2} \frac{\zeta^s}{(s+1)^2},
		\ee
		for all $r, t, s \geq 0$ where
		\bes
		\underline{C}_e = \overline{C}(h)(C_F + C_Q + C_R) > 0,
		\ees
		and $\overline{C}(h) > 0$ is a universal constant.
	\end{theorem}
	
	We now give the recursive estimate which is essential for our
	joint analyticity result.
	\begin{lemma}
		\label{Lemma:JointAnal:IndEst}
		For any integer $m > 0$, if $\tepsilon \in C^{\omega}_m(S_v)$
		such that 
		\bes
		\HolderNorm{\frac{\px^r \py^t \pz^s}{(r+t+s)!} \tepsilon }{m} 
		\leq C_{\tepsilon} \frac{\eta^r}{(r+1)^2} \frac{\theta^t}{(t+1)^2}
		\frac{\zeta^s}{(s+1)^2},
		\ees
		for all $r, t, s \geq 0$ and some constants 
		$C_{\tepsilon}, \eta, \theta, \zeta > 0$, and 
		\bes
		\SobNorm{ \frac{\px^r \py^t \pz^s}{(r+t+s)!} E_{\ell} }{2}
		\leq K B^{\ell} \frac{\eta^r} {(r+1)^2} \frac{\theta^t}{(t+1)^2}
		\frac{\zeta^s}{(s+1)^2},
		\quad \forall\ \ell < L,
		\ees
		for all $r, t, s \geq 0$ and for some constants $K, B > 0$. Then, 
		\bes
		\SobNorm{ \frac{\px^r \py^t \pz^s}{(r+t+s)!} F_L }{2} 
		\leq \tilde{C} K B^{L-1} \frac{\eta^r}{(r+1)^2}
		\frac{\theta^t}{(t+1)^2} \frac{\zeta^s}{(s+1)^2},
		\ees
		for all $r, t, s \geq 0$ and some constant $\tilde{C}>0$.
	\end{lemma}
	\begin{proof}
		Recall that
		\bes
		F_L = \bepsilon k_0^2 \tepsilon E_{L-1},
		\ees
		so, using Leibniz's rule we obtain,
		\bes
		\frac{\px^r \py^t \pz^s}{(r+t+s)!} F_L 
		= \frac{k^2_0 \bepsilon r! t! s!}{(r+t+s)!}
		\sum_{j=0}^{r} \sum_{k=0}^{t} \sum_{\ell=0}^{s}
		\left( \frac{\px^{r-j}}{(r-j)!} \frac{\py^{t-k}}{(t-k)!}
		\frac{\pz^{s-\ell}}{(s-\ell)!} \tepsilon \right)
		\left( \frac{\px^j}{j!} \frac{\py^k}{k!} \frac{\pz^{\ell}}{\ell!}
		E_{L-1} \right).
		\ees
		Using the inequality $r!t!s! \leq (r+t+s)!$ we obtain
		\begin{align*}
			\SobNorm{ \frac{\px^r \py^t \pz^s}{(r+t+s)!} F_L }{2}
			& \leq \bepsilon k_0^2 
			\sum_{j=0}^{r} \sum_{k=0}^{t} \sum_{\ell=0}^{s}
			\SobNorm{ \left( \frac{\px^{r-j}}{(r-j)!}
				\frac{\py^{t-k}}{(t-k)!} \frac{\pz^{s-\ell}}{(s-\ell)!}
				\tepsilon \right)
				\left( \frac{\px^j}{j!} \frac{\py^k}{k!} 
				\frac{\pz^\ell}{\ell!} E_{L-1} \right) }{2} \\
			& \leq \bepsilon k_0^2 \sum_{j=0}^{r} \sum_{k=0}^{t} \sum_{\ell=0}^{s}
			\tilde{M} \HolderNorm{ \frac{\px^{r-j}}{(r-j)!} \frac{\py^{t-k}}{(t-k)!}
				\frac{\pz^{s-\ell}}{(s-\ell)!} \tepsilon }{2}
			\\ & \quad \times
			\SobNorm{ \frac{\px^j}{j!} \frac{\py^k}{k!} \frac{\pz^\ell}{\ell!}
				E_{L-1} }{2} \\
			& \leq \bepsilon k_0^2 \tilde{M}
			\sum_{j=0}^{r} \sum_{k=0}^{t} \sum_{\ell=0}^{s} C_{\tepsilon} 
			\frac{\eta^{r-j}}{(r-j+1)^2} \frac{\theta^{t-k}}{(t-k+1)^2}
			\frac{\zeta^{s-\ell}}{(s-\ell+1)^2} \\
			& \quad \times K B^{L-1} \frac{\eta^j}{(j+1)^2}
			\frac{\theta^k}{(k+1)^2} \frac{\zeta^{\ell}}{(\ell+1)^2} \\
			& = \bepsilon k_0^2 C_{\tepsilon} \tilde{M} K B^{L-1} \frac{\eta^r}{(r+1)^2}
			\frac{\theta^t}{(t+1)^2} \frac{\zeta^s}{(s+1)^2} \\
			& \quad \times \sum_{j=0}^{r} \frac{(r+1)^2}{(r-j+1)^2 (j+1)^2}
			\sum_{k=0}^{t} \frac{(t+1)^2}{(t-k+1)^2(k+1)^2} \\
			& \quad \times \sum_{\ell=0}^{s} \frac{(s+1)^2}{(s-\ell+1)^2(\ell+1)^2} \\
			& \leq \bepsilon k_0^2 C_{\tepsilon} \tilde{M} S^3 K B^{L-1}
			\frac{\eta^r}{(r+1)^2} \frac{\theta^t}{(t+1)^2}
			\frac{\zeta^s}{(s+1)^2},
		\end{align*}
		where $S$ is a positive constant such that
		\bes
		\sum_{j=0}^{r} \frac{(r+1)^2}{(r-j+1)^2 (j+1)^2} < S,
		\quad \forall\ r \geq 0,
		\ees
		c.f.\ Lemma~\ref{Lemma:S}. Therefore, the proof is complete by choosing
		\bes
		\tilde{C} \geq \bepsilon k_0^2 C_{\tepsilon} \tilde{M} S^3.
		\ees
	\end{proof}
	
	We conclude with our joint analyticity theorem.
	\begin{theorem}
		Given any integer $m > 0$, if
		$(\omega,\bepsilon) \in \cP$, and
		$\tepsilon \in  C^{\omega}_m(S_v)$ such that 
		\bes
		\HolderNorm{ \frac{\px^r \py^t \pz^s}{(r+t+s)!} \tepsilon }{m}
		\leq C_{\tepsilon} \frac{\eta^r}{(r+1)^2} \frac{\theta^t}{(t+1)^2}
		\frac{\zeta^s}{(s+1)^2},
		\ees
		for all $r, t, s \geq 0$ and some constants 
		$C_{\tepsilon}, \eta, \theta, \zeta > 0$.
		Then the series \eqref{Eqn:EExp} converges strongly. Moreover the 
		$E_{\ell}(x,y,z)$ satisfy the joint analyticity estimate
		\be
		\label{Eqn:JointAnal}
		\SobNorm{ \frac{\px^r \py^t \pz^s}{(r+t+s)!} E_{\ell} }{2}
		\leq K B^{\ell} \frac{\eta^r}{(r+1)^2} \frac{\theta^t}{(t+1)^2}
		\frac{\zeta^s}{(s+1)^2},
		\ee
		for all $\ell, r, t, s \geq 0$ and some constants $K, B > 0$.
	\end{theorem}
	\begin{proof}
		We prove \eqref{Eqn:JointAnal} by induction, beginning with $\ell = 0$.
		Applying Theorem~\ref{Thm:JointAnal:EllEst} with $F \equiv 0$,
		$Q = \phi$, and $R \equiv 0$ we obtain
		\bes
		\SobNorm{ \frac{\px^r \py^t \pz^s}{(r+t+s)!} E_0 }{2} 
		\leq K \frac{\eta^r}{(r+1)^2} \frac{\theta^t}{(t+1)^2}
		\frac{\zeta^s}{(s+1)^2},
		\ees
		for all $r, t, s \geq 0$, where $K = \overline{C}(h) C_{\phi}$.
		Next we assume that \eqref{Eqn:JointAnal} is valid for all $\ell < L$.
		With $\ell = L$ we invoke Lemma~\ref{Lemma:JointAnal:IndEst} and 
		apply Theorem~\ref{Thm:JointAnal:EllEst} with $F \equiv F_L$,
		$C_F = \tilde{C} K B^{L-1}$, $Q \equiv 0$, and $R \equiv 0$,
		to arrive at
		\bes
		\SobNorm{ \frac{\px^r \py^t \pz^s}{(r+t+s)!} E_L }{2} 
		\leq \overline{C}(h) \tilde{C} K B^{L-1} \frac{\eta^r}{(r+1)^2}
		\frac{\theta^t}{(t+1)^2} \frac{\zeta^s}{(s+1)^2},
		\quad\ \forall r, t, s \geq 0.
		\ees
		The proof is complete by choosing $B > \overline{C}(h) \tilde{C}$.
	\end{proof}
	
	%
	%
	\appendix
	
	%
	%
	
	\section{Proof of the Elliptic Estimate}
	\label{Sec:EllEst:Proof}
	
	In this appendix we provide the proof of Theorem~\ref{Thm:EllEst} which
	has been so crucial to all of our developments.
	We will focus on establishing this result in the case $s=0$ as
	the case $s>0$ follows in analogous fashion.
	To begin, we recall that if the functions $v, F, Q, R$ from the
	theorem satisfy quasiperiodic
	boundary conditions, then they can each be expanded in 
	(generalized) Fourier series, e.g.,
	\bes
	v(x,y,z) = \sump \sumq \hat{v}_{p,q}(z) e^{i \alpha_p x + i \beta_q y},
	\quad
	\hat{v}_{p,q}(z) = \begin{pmatrix} \hat{v}^x_{p,q}(z) \\
		\hat{v}^y_{p,q}(z) \\ \hat{v}^z_{p,q}(z) \end{pmatrix}.
	\ees
	In terms of these expansions we have the following restatement
	of the governing equations \eqref{Eqn:EllEst}.
	\begin{lemma}
		\label{Lemma:BVP:vx:vy}
		Let $v$ be the solution of \eqref{Eqn:EllEst}. Under the assumptions of 
		Theorem~\ref{Thm:EllEst}, $\hat{v}^x_{p,q}$ and $\hat{v}^y_{p,q}$ satisfy
		the following systems of two point boundary value problems
		\bse
		\begin{align}
			& \pz^2 \hat{v}^x_{p,q} + (\gammae_{p,q})^2 \hat{v}^x_{p,q} = H_{p,q}(z),
			&& -h < z < h, 
			\label{Eqn:BVP:vx} \\
			& \pz \hat{v}^x_{p,q}(h) - i \gammau_{p,q} \hat{v}^x_{p,q}(h) = -\hat{Q}^x_{p,q}, \\
			& \pz \hat{v}^x_{p,q}(-h) + i \gammaw_{p,q} \hat{v}^x_{p,q}(-h) = \hat{R}^x_{p,q},
		\end{align}
		\ese
		and
		\bse
		\begin{align}
			& \pz^2 \hat{v}^y_{p,q} + (\gammae_{p,q})^2 \hat{v}^y_{p,q} = L_{p,q}(z),
			&& -h < z < h,
			\label{Eqn:BVP:vy} \\
			& \pz \hat{v}^y_{p,q}(h) - i \gammau_{p,q} \hat{v}^y_{p,q}(h) = -\hat{Q}^y_{p,q}, \\
			& \pz \hat{v}^y_{p,q}(-h) + i \gammaw_{p,q} \hat{v}^y_{p,q}(-h) = \hat{R}^y_{p,q},
		\end{align}
		\ese
		where
		\bes
		\gammae_{p,q} := \begin{cases}
			\sqrt{\bepsilon k_0^2 - \alpha_p^2 - \beta_q^2},
			& \alpha_p^2 + \beta_q^2 \leq \bepsilon k_0^2, \\
			i \sqrt{\alpha_p^2 + \beta_q^2 - \bepsilon k_0^2}
			:= i \bgamma_{p,q},
			& \alpha_p^2 + \beta_q^2 > \bepsilon k_0^2,
		\end{cases}
		\ees
		and
		\bes
		\gammam_{p,q} := \begin{cases}
			\sqrt{\epsilon^{(m)} k_0^2 - \alpha_p^2 - \beta_q^2},
			& \alpha_p^2 + \beta_q^2 \leq \epsilon^{(m)} k_0^2, \\
			i \sqrt{\alpha_p^2 + \beta_q^2 - \epsilon^{(m)} k_0^2}
			:= i \bgamma_{p,q}^{(m)},
			& \alpha_p^2 + \beta_q^2 > \epsilon^{(m)} k_0^2,
		\end{cases}
		\ees
		for $m \in \{u , v\}$, and $\bgamma_{p,q}, \bgamma_{p,q}^{(m)} \in \Real^+$ and
		\bse
		\begin{align}
			H_{p,q}(z) & := \frac{\alpha_p \beta_q}{ k_0^2 \bepsilon } \hat{F}^y_{p,q} 
			+ \frac{\alpha_p^2 - k_0^2 \bepsilon}{ k_0^2 \bepsilon } \hat{F}^x_{p,q}
			- \frac{i \alpha_p}{ k_0^2 \bepsilon} \pz \hat{F}^z_{p,q},
			\label{Eqn:H} \\
			L_{p,q}(z) & := \frac{\alpha_p \beta_q}{ k_0^2 \bepsilon } \hat{F}^x_{p,q}
			+ \frac{\beta_q^2 - k_0^2 \bepsilon}{ k_0^2 \bepsilon } \hat{F}^y_{p,q}
			- \frac{i \beta_{p,q}}{ k_0^2 \bepsilon} \pz \hat{F}^z_{p,q}.
			\label{Eqn:L}
		\end{align}
		\ese
		Furthermore, we can compute $\hat{v}^z_{p,q}$ from these as
		\be
		\label{Eqn:vz}
		\hat{v}^z_{p,q} = - \frac{1}{(\gammae_{p,q})^2} \hat{F}^z_{p,q} 
		+ \frac{i \alpha_p}{(\gammae_{p,q})^2} \pz \hat{v}^x_{p,q} 
		+ \frac{i \beta_q}{(\gammae_{p,q})^2} \pz \hat{v}^y_{p,q}.
		\ee
	\end{lemma}
	\begin{proof}
		We begin with the observation that
		\bes
		\Curl{\Curl{v}} = -\Laplacian{v} + \Grad{\Div{v}},
		\ees
		so that
		\bes
		\cL_0[v] = -\Laplacian{v} + \Grad{\Div{v}} 
		- k_0^2 \bepsilon.
		\ees
		Next we apply $\cL_0$ to the expansion
		\bes
		\cL_0[v] = \sump \sumq \cL_0[\hat{v}_{p,q}(z)
		e^{i \alpha_p x + i \beta_q y}],
		\ees
		which requires
		\bes
		-\Laplacian{\hat{v}^m_{p,q} 
			e^{i \alpha_p x + i \beta_q y} } 
		= \left\{ (\alpha_p^2 + \beta_q^2) \hat{v}^m_{p,q}
		- \pz^2 \hat{v}^m_{p,q} \right\} 
		e^{i \alpha_p x + i \beta_q y},
		\quad
		m \in \{ x, y, z \},
		\ees
		and
		\bes
		\Div{ \hat{v}_{p,q} e^{i \alpha_p x + i \beta_q y} }
		= \left\{ (i \alpha_p) \hat{v}^x_{p,q} 
		+ (i \beta_q) \hat{v}^y_{p,q} 
		+ \pz \hat{v}^z_{p,q} \right\} 
		e^{i \alpha_p x + i \beta_q y},
		\ees
		and
		\begin{align*}
			\px \Div{\hat{v}_{p,q} e^{i \alpha_p x + i \beta_q y}}
			& = \left\{ -\alpha_p^2 \hat{v}^x_{p,q} 
			- \alpha_p \beta_q \hat{v}^y_{p,q} 
			+ (i \alpha_p) \pz \hat{v}^z_{p,q} \right\} 
			e^{i \alpha_p x + i \beta_q y}, \\
			\py \Div{\hat{v}_{p,q} e^{i \alpha_p x + i \beta_q y}}
			& = \left\{ -\alpha_p \beta_q \hat{v}^x_{p,q} 
			- \beta_q^2 \hat{v}^y_{p,q} 
			+ (i \beta_q) \pz \hat{v}^z_{p,q} \right\} 
			e^{i \alpha_p x + i \beta_q y}, \\
			\pz \Div{\hat{v}_{p,q} e^{i \alpha_p x + i \beta_q y}}
			& = \left\{ (i \alpha_p) \pz \hat{v}^x_{p,q} 
			+ (i \beta_q) \pz \hat{v}^y_{p,q} 
			+ \pz^2 \hat{v}^z_{p,q} \right\} 
			e^{i \alpha_p x + i \beta_q y}.
		\end{align*}
		From these we have
		\bse
		\label{Eqn:ODEs}
		\begin{align}
			& -\alpha_p \beta_q \hat{v}^y_{p,q}
			+ \beta_q^2 \hat{v}^x_{p,q}
			- \pz^2 \hat{v}^x_{p,q}
			+ (i \alpha_p) \pz \hat{v}^z_{p,q}
			- k_0^2 \bepsilon \hat{v}^x_{p,q}
			= \hat{F}^x_{p,q} 
			\label{Eqn:ODEs:a} \\
			& -\alpha_p \beta_q \hat{v}^x_{p,q}
			+ \alpha_p^2 \hat{v}^y_{p,q}
			- \pz^2 \hat{v}^y_{p,q}
			+ (i \beta_q) \pz \hat{v}^z_{p,q}
			- k_0^2 \bepsilon \hat{v}^y_{p,q}
			= \hat{F}^y_{p,q}
			\label{Eqn:ODEs:b} \\
			& -(\gammae_{p,q})^2 \hat{v}^z_{p,q}
			+ (i \alpha_p) \pz \hat{v}^x_{p,q}
			+ (i \beta_q) \pz \hat{v}^y_{p,q}
			= \hat{F}^z_{p,q}.
			\label{Eqn:ODEs:c}
		\end{align}
		\ese
		If we now multiply \eqref{Eqn:ODEs:a} by $\beta_q$ 
		and \eqref{Eqn:ODEs:b} by $\alpha_p$, and subtract
		them we obtain
		\bes
		-\beta_q (\pz^2 \hat{v}^x_{p,q} 
		+ (\gammae_{p,q})^2 \hat{v}^x_{p,q} ) 
		+ \alpha_p (\pz^2 \hat{v}^y_{p,q} 
		+ (\gammae_{p,q})^2 \hat{v}^y_{p,q})
		= \beta_q \hat{F}^x_{p,q} 
		- \alpha_p \hat{F}^y_{p,q}.
		\ees
		Furthermore, dividing \eqref{Eqn:ODEs:c} by 
		$(\gammae_{p,q})^2$ and then differentiating the result
		with respect to $z$, we obtain
		\bes
		\pz \hat{v}^z_{p,q} = \frac{1}{(\gammae_{p,q})^2}
		\left( (i \alpha_p) \pz^2 \hat{v}^x_{p,q} 
		+ (i \beta_q) \pz^2 \hat{v}^y_{p,q} - \pz \hat{F}^z_{p,q}\right).
		\ees
		Substituting this into \eqref{Eqn:ODEs:b} we obtain
		\bes
		\alpha_p \beta_q (\pz^2 \hat{v}^x_{p,q} 
		+ (\gammae_{p,q})^2 \hat{v}^x_{p,q}) 
		- (\alpha_p^2 - k_0^2 \bepsilon) (\pz^2 \hat{v}^y_{p,q}
		+ (\gammae_{p,q})^2 \hat{v}^y_{p,q}) 
		= - (i \beta_q) \pz \hat{F}^z_{p,q}-(\gammae_{p,q})^2 \hat{F}^y_{p,q}.
		\ees
		If we denote
		\bes
		U := \pz^2 \hat{v}^x_{p,q} + (\gammae_{p,q})^2 \hat{v}^x_{p,q},
		\quad
		W := \pz^2 \hat{v}^y_{p,q} + (\gammae_{p,q})^2 \hat{v}^y_{2n},
		\ees
		then we find a system of equations for $U$ and $W$
		\begin{align*}
			-\beta_q U + \alpha_p W & = \beta_q \hat{F}^x_{p,q} 
			- \alpha_p \hat{F}^y_{p,q}, \\
			\alpha_p \beta_q U - (\alpha_p^2 - k_0^2 \bepsilon) W 
			& = - (i \beta_q) \pz \hat{F}^z_{p,q}-(\gammae_{p,q})^2 \hat{F}^y_{p,q}.
		\end{align*}
		Solving this system gives us
		\begin{align*}
			U & = \frac{\alpha_p \beta_q}{k_0^2 \bepsilon} \hat{F}^y_{p,q}
			+ \frac{\alpha_p^2 - k_0^2 \bepsilon}{k_0^2 \bepsilon} \hat{F}^x_{p,q} 
			- \frac{i \alpha_p}{k_0^2 \bepsilon} \pz \hat{F}^z_{p,q} =: H_{p,q}(z),
			\\
			W & = \frac{\alpha_p \beta_q}{k_0^2 \bepsilon} \hat{F}^x_{p,q}
			+ \frac{\beta_q^2 - k_0^2 \bepsilon}{k_0^2 \bepsilon} \hat{F}^y_{p,q}
			- \frac{i \beta_q}{k_0^2 \bepsilon} \pz \hat{F}^z_{p,q} =: L_{p,q}(z). 
		\end{align*}
		The proof is complete.
	\end{proof}
	
	\begin{lemma}
		The solutions $\hat{v}^x_{p,q}$ and $\hat{v}^y_{p,q}$ of \eqref{Eqn:ODEs:a}
		and \eqref{Eqn:ODEs:b} are uniquely determined by
		\begin{align*}
			\hat{v}^x_{p,q}(z) & = -\hat{Q}^x_{p,q} \phi_{h}(z;p,q) 
			- \hat{R}^x_{p,q} \phi_{-h}(z;p,q) \frac{2 \bgamma_{p,q}}{D} 
			- I_h[H_{p,q}](z) - I_{-h}[H_{p,q}](z), \\
			\hat{v}^y_{p,q}(z) & = -\hat{Q}^y_{p,q} \phi_{h}(z;p,q) 
			- \hat{R}^y_{p,q} \phi_{-h}(z;p,q) \frac{2 \bgamma_{p,q}}{D}
			- I_h[L_{p,q}](z) - I_{-h}[L_{p,q}](z),
		\end{align*}
		where 
		\begin{align*}
			\phi_{h}(z;p,q) & := \left( \frac{\bgamma_{p,q} + \bgamma_{p,q}^{(w)}}{D} 
			\right) e^{\bgamma_{p,q}(z+h)}
			+ \left( \frac{\bgamma_{p,q} - \bgamma_{p,q}^{(w)}}{D} \right)
			e^{-\bgamma_{p,q}(z+h)}, \\
			\phi_{-h}(z;p,q) & := \left( \frac{\bgamma_{p,q} - \bgamma_{p,q}^{(u)}}
			{2 \bgamma_{p,q}} \right) e^{\bgamma_{p,q}(z-h)} 
			+ \left( \frac{\bgamma_{p,q} + \bgamma_{p,q}^{(u)}}{2 \bgamma_{p,q}} \right)
			e^{-\bgamma_{p,q}(z-h)}, \\
			I_h[\zeta](z) & := \int_z^h \phi_{h}(z;p,q) \phi_{-h}(s;p,q) \zeta(s) \; ds, \\
			I_{-h}[\zeta](z) & := \int_{-h}^z \phi_{-h}(z;p,q) \phi_{h}(s;p,q) \zeta(s) \; ds, \\
			D & := (\bgamma_{p,q} + \bgamma_{p,q}^{(w)}) 
			(\bgamma_{p,q} + \bgamma_{p,q}^{(u)}) e^{2 \bgamma_{p,q} h}	
			- (\bgamma_{p,q} - \bgamma_{p,q}^{(u)}) 
			(\bgamma_{p,q} - \bgamma_{p,q}^{(w)}) e^{-2 \bgamma_{p,q} h}.
		\end{align*}
	\end{lemma}
	
	With this we are ready to give the proof of Theorem~\ref{Thm:EllEst}. As stated
	above, we provide a detailed proof for the estimate \eqref{Eqn:EllEst:Est}
	in the case when $s=0$.
	
	\begin{proof}{[Theorem~\ref{Thm:EllEst}]}  
		To start we recall that
		\begin{multline*}
			\SobNorm{v}{2}^2 :=
			\sump \sumq \left( \Angle{(p,q)}^4 \int_{-h}^h \Abs{\hat{v}_{p,q}(z)}^2 \; dz
			+ \Angle{(p,q)}^2 \int_{-h}^h \Abs{\pz \hat{v}_{p,q}(z)}^2 \; dz 
			\right. \\ \left.
			+ \int_{-h}^h \Abs{\pz^2 \hat{v}_{p,q}(z)}^2 \; dz \right).
		\end{multline*}
		We point out that the indices in the double sum on $(p,q)$ 
		can be divided into
		two sets: The \textit{propagating modes} which are defined by
		\bes
		\mathbf{\bar{P}} := \left\{ (p,q) \in \Integer^2\ |\ 
		\alpha_p^2 + \beta_q^2 \leq \bepsilon k_0^2 \right\},
		\ees
		and the \textit{evanescent modes} specified by
		\bes
		\mathbf{\bar{E}} := \left\{ (p,q) \in \Integer^2\ |\ 
		\alpha_p^2 + \beta_q^2 > \bepsilon k_0^2 \right\}.
		\ees
		The former is of \textit{finite} size and gives \textit{complex}
		$\bgamma_{p,q}$, while the latter is \textit{unbounded} and
		features $\bgamma_{p,q}$ real and positive. From the previous
		Lemma we observe that
		$\{ \hat{v}_{p,q}(z), \pz \hat{v}_{p,q}(z), \pz^2 \hat{v}_{p,q}(z) \}$ 
		are all bounded on $-h < z < h$ so that there exists a constant
		$K_0 > 0$ such that, among the propagating modes,
		\bes
		\max \left\{ \Norm{\hat{v}_{p,q}}{L^2}, \Norm{\pz \hat{v}_{p,q}}{L^2},
		\Norm{\pz^2 \hat{v}_{p,q}}{L^2} \right\} < K_0,
		\quad
		\forall\ (p,q) \in \mathbf{\bar{P}}.
		\ees
		As there are only a finite number of these, we can estimate all of them
		uniformly by 
		\bes
		\max_{(p,q) \in \mathbf{\bar{P}}} \left\{
		\Abs{\hat{Q}_{p,q}}, \Abs{\hat{R}_{p,q}},
		\Norm{\hat{H}_{p,q}}{L^2}, \Norm{\hat{L}_{p,q}}{L^2}
		\right\}.
		\ees
		For this reason, we restrict our subsequent developments to the
		evanescent modes, which requires a careful asymptotic study
		as $\Abs{(p,q)}$ grows, and assume that $\bgamma_{p,q}$ are real
		and positive.
		
		We begin by estimating $\hat{v}^x_{p,q}$. If we denote
		\bes
		a_{p,q} = \bgamma_{p,q} + \bgamma_{p,q}^{(w)},
		\quad
		b_{p,q} = \bgamma_{p,q} - \bgamma_{p,q}^{(w)},
		\quad
		c_{p,q} = \bgamma_{p,q} - \bgamma_{p,q}^{(u)},
		\quad
		d_{p,q} = \bgamma_{p,q} + \bgamma_{p,q}^{(u)},
		\ees
		then,
		\begin{align*}
			\phi_{h}(z;p,q) & = \frac{a_{p,q}}{D} e^{\bgamma_{p,q}(z+h)}
			+ \frac{b_{p,q}}{D} e^{-\bgamma_{p,q}(z+h)}, \\
			\phi_{-h}(z;p,q) & = \frac{c_{p,q}}{2 \bgamma_{p,q}} e^{\bgamma_{p,q}(z-h)} 
			+ \frac{d_{p,q}}{2 \bgamma_{p,q}} 
			e^{-\bgamma_{p,q}(z-h)}, \\
			D & = a_{p,q} d_{p,q} e^{2 \bgamma_{p,q} h}	
			- b_{p,q} c_{p,q} e^{-2 \bgamma_{p,q} h},
		\end{align*}
		and $\hat{v}^x_{p,q}(z) = \sum_{j=1}^{10} S_j(z)$, where
		\begin{gather*}
			S_1(z) = -\hat{Q}^x_{p,q} \frac{a_{p,q}}{D} e^{\bgamma_{p,q}(z+h)},
			\quad
			S_2(z) = -\hat{Q}^x_{p,q} \frac{b_{p,q}}{D} e^{-\bgamma_{p,q}(z+h)},
			\\
			S_3(z) = -\hat{R}^x_{p,q} \frac{c_{p,q}}{2 \bgamma_{p,q}} 
			e^{\bgamma_{p,q}(z-h)},
			\quad
			S_4(z) = -\hat{R}^x_{p,q} \frac{d_{p,q}}{2 \bgamma_{p,q}}
			e^{-\bgamma_{p,q}(z-h)}, \\
			S_5(z) = -\frac{a_{p,q} c_{p,q}}{2 \bgamma_{p,q} D}
			\int_{-h}^h e^{\bgamma_{p,q}(z+s)} H_{p,q}(s) \; ds \\
			S_6(z) = -\frac{b_{p,q} d_{p,q}}{2 \bgamma_{p,q} D}
			\int_{-h}^h e^{-\bgamma_{p,q}(z+s)} H_{p,q}(s) \; ds \\
			S_7(z) = -\frac{a_{p,q} d_{p,q}}{2 \bgamma_{p,q} D}
			\int_{z}^h e^{\bgamma_{p,q}(z-s+2h)} H_{p,q}(s) \; ds \\
			S_8(z) = -\frac{b_{p,q} c_{p,q}}{2 \bgamma_{p,q} D}
			\int_{-h}^z e^{\bgamma_{p,q}(z-s-2h)} H_{p,q}(s) \; ds \\
			S_9(z) = -\frac{b_{p,q} c_{p,q}}{2 \bgamma_{p,q} D}
			\int_{z}^h e^{-\bgamma_{p,q} (z-s+2h)} H_{p,q}(s) \; ds \\
			S_{10}(z) = - \frac{a_{p,q} d_{p,q}}{2 \bgamma_{p,q} D}
			\int_{-h}^z e^{-\bgamma_{p,q}(z-s-2h)} H_{p,q}(s) \; ds.
		\end{gather*}
		To estimate $\Norm{\hat{v}^x_{p,q}}{L^2}$ one must address each of
		these ten terms individually and use
		\be
		\label{Eqn:T1:T10}
		\Norm{\hat{v}^x_{p,q}}{L^2} 
		\leq \sum_{j=1}^{10} T_j,
		\quad
		T_j := \Norm{S_j}{L^2}.
		\ee
		For brevity we provide details on two of these, $T_1$
		and $T_5$. For the former we begin
		\begin{align*}
			T_1 & \leq \Abs{\hat{Q}^x_{p,q}} 
			\left( \Abs{ \frac{a_{p,q}}{D} }^2 
			\int_{-h}^h e^{2 \bgamma_{p,q}(z+h)} \; dz \right)^{1/2} \\
			& = \Abs{\hat{Q}^x_{p,q}} 
			\left( \Abs{ \frac{a_{p,q}}{D} }^2
			\frac{e^{4 \bgamma_{p,q} h} - 1}{2 \bgamma_{p,q}} \right)^{1/2} \\
			& = \frac{\Abs{\hat{Q}^x_{p,q}}}{\sqrt{2 \bgamma_{p,q}} \Abs{d_{p,q}}}
			\left( \frac{(e^{4 \bgamma_{p,q} h} - 1)^{1/2}}
			{\Abs{e^{2 \bgamma_{p,q} h} - \frac{b_{p,q} c_{p,q}}{a_{p,q} d_{p,q}}
					e^{-2 \bgamma_{p,q} h}}} \right).
		\end{align*}
		Since
		\bes
		\lim \limits_{\Abs{(p,q)} \rightarrow \infty} 
		\frac{(e^{4 \bgamma_{p,q} h} - 1)^{1/2}}{\Abs{e^{2 \bgamma_{p,q} h}
				- \frac{b_{p,q} c_{p,q}}{a_{p,q} d_{p,q}} e^{-2 \bgamma_{p,q} h}}} 
		= \lim \limits_{\Abs{(p,q)} \rightarrow \infty}
		\frac{(1- e^{-4 \bgamma_{p,q} h})^{1/2}}
		{\Abs{1-\frac{b_{p,q} c_{p,q}}{a_{p,q} d_{p,q}} 
				e^{-4 \bgamma_{p,q} h}}} = 1,
		\ees
		there exists a constant $C > 0$ such that 
		\bes
		\frac{(e^{4 \bgamma_{p,q} h} - 1)^{1/2}}
		{\Abs{e^{2 \bgamma_{p,q} h} - \frac{b_{p,q} c_{p,q}}
				{a_{p,q} d_{p,q}} e^{-2 \bgamma_{p,q} h}}} 
		\leq C 
		\quad \forall\ (p,q) \in \Integer^2,
		\ees
		therefore,
		\bes
		T_1 \leq C \frac{\Abs{\hat{Q}^x_{p,q}}}
		{\sqrt{2 \bgamma_{p,q}} \Abs{d_{p,q}}}.
		\ees
		For the latter, by using H\"older's inequality, we obtain
		\begin{align*}
			T_5 & = \Abs{ \frac{a_{p,q} c_{p,q}}{2 \bgamma_{p,q} D} }
			\left( \int_{-h}^h \Abs{ \int_{-h}^h e^{\bgamma_{p,q} (z+s)} 
				H_{p,q}(s) \; ds }^2 \; dz \right)^{1/2} \\
			& \leq \Abs{ \frac{a_{p,q} c_{p,q}}{2 \bgamma_{p,q} D} }
			\left( \int_{-h}^h \left( \int_{-h}^h e^{2 \bgamma_{p,q}(z+s)} \; ds \right) 
			\left( \int_{-h}^h \Abs{H_{p,q}(s)}^2 \; ds \right) \; dz \right)^{1/2} \\
			& = \Abs{ \frac{a_{p,q} c_{p,q}}{2 \bgamma_{p,q} D} }
			\left( \int_{-h}^{h} \frac{e^{2 \bgamma_{p,q} (z+h)} 
				- e^{2 \bgamma_{p,q}(z-h)}}{2 \bgamma_{p,q}} \; dz \right)^{1/2}
			\Norm{H_{p,q}}{L^2} \\
			& = \Abs{ \frac{a_{p,q} c_{p,q}}{2 \bgamma_{p,q} D} }
			\left( \frac{e^{4 \bgamma_{p,q} h} + e^{-4 \bgamma_{p,q} h} - 2}
			{4 \bgamma_{p,q} ^2} \right)^{1/2} \Norm{H_{p,q}}{L^2} \\
			& = \frac{( e^{4 \bgamma_{p,q} h}
				+ e^{-4 \bgamma_{p,q} h} - 2 )^{1/2}}
			{\Abs{e^{2 \bgamma_{p,q} h} - \frac{b_{p,q} c_{p,q}}
					{a_{p,q} d_{p,q}} \; e^{-2 \bgamma_{p,q} h}}} \;
			\frac{\Abs{c_{p,q}} \Norm{H_{p,q}}{L^2}}{4 \bgamma_{p,q}^2 \Abs{d_{p,q}}} \\
			& = \frac{(1 + e^{-8 \bgamma_{p,q} h} - 2 e^{-4 \bgamma_{p,q} h})^{1/2}}
			{\Abs{1 - \frac{b_{p,q} c_{p,q}}{a_{p,q} d_{p,q}} \;
					e^{-4 \bgamma_{p,q} h}}} \; \frac{\Abs{c_{p,q}} \Norm{H_{p,q}}{L^2}}
			{4 \bgamma_{p,q}^2 \Abs{d_{p,q}}}
			\leq C \frac{\Norm{H_{p,q}}{L^2}}{4 \bgamma_{p,q}^2 \Abs{d_{p,q}}},
		\end{align*}
		where the last inequality was obtained by using the fact that 
		\bes
		\lim \limits_{\Abs{(p,q) \rightarrow \infty}}
		\frac{(1 + e^{-8 \bgamma_{p,q} h} - 2 e^{-4 \bgamma_{p,q} h})^{1/2}}
		{\Abs{1 - \frac{b_{p,q} c_{p,q}}{a_{p,q} d_{p,q}} \; e^{-4 \bgamma_{p,q} h}}}
		= 1.
		\ees

		Substituting the estimates of $T_1, \ldots, T_{10}$ into
		\eqref{Eqn:T1:T10} we obtain
		\begin{align}
			\label{Eqn:vx:Est}
			\Norm{\hat{v}^x_{p,q}}{L^2} &\leq 
			C \frac{\Abs{\hat{Q}^x_{p,q}}}{\sqrt{\bgamma_{p,q}} \Abs{d_{p,q}}}
			+ C \frac{\Abs{\hat{R}^x_{p,q}}}{\sqrt{\bgamma_{p,q}} \Abs{a_{p,q}}}
			+ C \frac{\Norm{H_{p,q}}{L^2}}{\bgamma_{p,q}^2 \Abs{d_{p,q}}} 
			+ C \frac{\Norm{H_{p,q}}{L^2}}{\bgamma_{p,q}^2} \\ \nonumber
			& \leq 
			C \frac{\Abs{\hat{Q}^x_{p,q}}}{\sqrt{\bgamma_{p,q}} \Abs{d_{p,q}}}
			+ C \frac{\Abs{\hat{R}^x_{p,q}}}{\sqrt{\bgamma_{p,q}} \Abs{a_{p,q}}}
			+ C \frac{\Norm{H_{p,q}}{L^2}}{\bgamma_{p,q}^2}.
		\end{align}
		Similarly, we obtain
		\be
		\label{Eqn:vy:Est}
		\Norm{\hat{v}^y_{p,q}}{L^2} \leq 
		C \frac{\Abs{\hat{Q}^y_{p,q}}}{\sqrt{\bgamma_{p,q}} \Abs{d_{p,q}}}
		+ C \frac{\Abs{\hat{R}^y_{p,q}}}{\sqrt{\bgamma_{p,q}} \Abs{a_{p,q}}} 
		+ C \frac{\Norm{L_{p,q}}{L^2}}{\bgamma_{p,q}^2 }.
		\ee
		
		Next, we estimate $\hat{v}^z_{p,q}$ and from \eqref{Eqn:vz} we obtain
		\begin{align*}
			\hat{v}^z_{p,q}
			& = \frac{1}{\bgamma_{p,q}^2} \hat{F}^z_{p,q} 
			- \frac{i \alpha_p}{\bgamma_{p,q}^2} \pz \hat{v}^x_{p,q} 
			- \frac{i \beta_q}{\bgamma_{p,q}^2} \pz \hat{v}^y_{p,q} \\
			& = \frac{1}{\bgamma_{p,q}^2} \hat{F}^z_{p,q} 
			- \frac{i \alpha_p}{\bgamma_{p,q}^2}
			\left[ -\hat{Q}^x_{p,q} \pz \phi_{h} 
			- \hat{R}^x_{p,q} \pz \phi_{-h}
			\frac{2 \bgamma_{p,q}}{D} - \pz I_h[H_{p,q}] - \pz I_{-h}[H_{p,q}] \right] \\
			& \quad
			- \frac{i \beta_q}{\bgamma_{p,q}^2}
			\left[ -\hat{Q}^y_{p,q} \pz \phi_{h} 
			- \hat{R}^y_{p,q} \pz \phi_{-h}
			\frac{2 \bgamma_{p,q}}{D} - \pz I_h[L_{p,q}] - \pz I_{-h}[L_{p,q}] \right] \\
			& = \frac{1}{\bgamma_{p,q}^2} \hat{F}^z_{p,q}
			+ \frac{i \pz \phi_{h}}{\bgamma_{p,q}^2}
			\left( \alpha_p \hat{Q}^x_{p,q} + \beta_q \hat{Q}^y_{p,q} \right)
			+ \frac{2 i \pz \phi_{-h}}{\bgamma_{p,q} D}
			\left( \alpha_p \hat{R}^x_{p,q} 
			+ \beta_q \hat{R}^y_{p,q} \right) \\
			& \quad
			+ \frac{i \alpha_p}{\bgamma_{p,q}^2}
			\left( \pz I_h[H_{p,q}] + \pz I_{-h}[H_{p,q}] \right)
			+ \frac{i \beta_q}{\bgamma_{p,q}^2}
			\left( \pz I_h[L_{p,q}] + \pz I_{-h}[L_{p,q}] \right).
		\end{align*}
		Therefore,
		\bes
		\Norm{\hat{v}^z_{p,q}}{L^2} 
		= \left( \int_{-h}^h \Abs{\hat{v}^z_{p,q}(z)}^2 \; dz 
		\right)^{1/2}
		\leq \frac{1}{\bgamma_{p,q}^2} \Norm{\hat{F}^z_{p,q}}{L^2}
		+ \sum_{j=1}^{4} I_j,
		\ees
		where
		\begin{align*}
			I_1 & :=  \frac{\Abs{\alpha_p \hat{Q}^x_{p,q} + \beta_q \hat{Q}^y_{p,q}}}
			{\bgamma_{p,q}^2}
			\left( \int_{-h}^h \Abs{\pz \phi_{h}}^2 \; dz \right)^{1/2}, \\
			I_2 & := \frac{2 \Abs{\alpha_p \hat{R}^x_{p,q} + \beta_q \hat{R}^y_{p,q}}}
			{\bgamma_{p,q} \Abs{D}}
			\left( \int_{-h}^h \Abs{\pz \phi_{-h}}^2 \; dz \right)^{1/2}, \\
			I_3 & := \frac{\Abs{\alpha_p}}{\bgamma_{p,q}^2}
			\left( \int_{-h}^h \Abs{ \pz I_h[H_{p,q}] + \pz I_{-h}[H_{p,q}] }^2 \;
			dz \right)^{1/2}, \\
			I_4 & := \frac{\Abs{\beta_q}}{\bgamma_{p,q}^2}
			\left( \int_{-h}^h \Abs{\pz I_h[L_{p,q}] + \pz I_{-h}[L_{p,q}]}^2 \;
			dz \right)^{1/2}.
		\end{align*}
		All four of these $I_j$ must be estimated, but we focus on $I_1$ and
		$I_3$ to streamline our presentation. To start,
		\bes
		\pz \phi_{h}(z;p,q) = \frac{a_{p,q} \bgamma_{p,q}}{D}
		e^{\bgamma_{p,q}(z+h)} - \frac{b_{p,q} \bgamma_{p,q}}{D}
		e^{-\bgamma_{p,q}(z+h)},
		\ees
		so, by the H\"older Inequality, we have, cancelling a factor
		of $\bgamma_{p,q}$,
		\bes
		I_1 \leq \frac{\Abs{\alpha_p \hat{Q}^x_{p,q} + \beta_q \hat{Q}^y_{p,q}}}
		{\bgamma_{p,q}} \left( 
		\Abs{\frac{a_{p,q}}{D}} \Norm{e^{\bgamma (z+h)}}{L^2} + 
		\Abs{\frac{b_{p,q}}{D}} \Norm{e^{-\bgamma (z+h)}}{L^2} \right).
		\ees
		The discrete Cauchy--Schwartz inequality tells us that
		\bes
		\frac{\Abs{\alpha_p \hat{Q}^x_{p,q} + \beta_q \hat{Q}^y_{p,q}}}
		{\bgamma_{p,q}}
		\leq \frac{\sqrt{\alpha_p^2 + \beta_q^2} \Abs{\hat{Q}_{p,q}}}
		{\bgamma_{p,q}}
		\leq C \Abs{\hat{Q}_{p,q}},
		\ees
		where the final inequality comes from
		\bes
		\lim \limits_{\Abs{(p,q)} \rightarrow \infty}
		\frac{\sqrt{\alpha_p^2 + \beta_q^2}}{\bgamma_{p,q}} = 1.
		\ees
		So, we continue,
		\begin{align*}
			I_1 & \leq C \Abs{\hat{Q}_{p,q}} \left\{ 
			\Abs{\frac{a_{p,q}}{D}} \left( \frac{e^{4 \bgamma_{p,q} h - 1}}
			{2 \bgamma_{p,q}} \right)^{1/2}
			+ \Abs{\frac{b_{p,q}}{D}} \left( \frac{1 - e^{-4 \bgamma_{p,q}}}
			{2 \bgamma_{p,q}} \right)^{1/2} \right\} \\
			& \leq C \frac{\Abs{\hat{Q}_{p,q}}}{\sqrt{2 \bgamma_{p,q}}}
			\left\{ 
			\frac{1}{\Abs{d_{p,q}}}
			\frac{(1 - e^{-4 \bgamma_{p,q} h})^{1/2}}
			{\Abs{1 - \frac{b_{p,q} c_{p,q}}{a_{p,q} d_{p,q}}
					e^{-4 \bgamma_{p,q} h}}}
			\right. \\ & \quad \left.
			+ \frac{\Abs{b_{p,q}}}{\Abs{a_{p,q} d_{p,q}} 
				e^{2 \bgamma_{p,q} h}}
			\frac{(1 - e^{-4 \bgamma_{p,q} h})^{1/2}}
			{\Abs{1 - \frac{b_{p,q} c_{p,q}}{a_{p,q} d_{p,q}}
					e^{-4 \bgamma_{p,q} h}}}
			\right\}.
		\end{align*}
		Since
		\bes
		\lim \limits_{\Abs{(p,q)} \rightarrow \infty}
		\frac{(1 - e^{-4 \bgamma_{p,q} h})^{1/2}}
		{\Abs{1 - \frac{b_{p,q} c_{p,q}}{a_{p,q} d_{p,q}}
				e^{-4 \bgamma_{p,q} h}}} = 1,
		\ees
		we have that
		\bes
		I_1 \leq C \frac{\Abs{\hat{Q}_{p,q}}}
		{\sqrt{2 \bgamma_{p,q}} \Abs{d_{p,q}}}.
		\ees
		Continuing, we have
		\begin{align*}
			\pz I_h[H_{p,q}](z) & = \pz \phi_{h}(z) \int_z^h \phi_{-h}(s) H_{p,q}(s) \; ds
			- \phi_{h}(z) \phi_{-h}(z) H_{p,q}(z), \\
			\pz I_{-h}[H_{p,q}](z) & = \pz \phi_{-h}(z) \int_{-h}^z \phi_{h}(s) H_{p,q}(s) \; ds
			+ \phi_{-h}(z) \phi_{h}(z) H_{p,q}(z).	
		\end{align*}
		so that
		\begin{multline}
			\label{Eqn:pzIh}
			\pz I_h[H_{p,q}](z) + \p_z I_{-h}[H_{p,q}](z)
			\\
			= \int_z^h \pz \phi_{h}(z) \phi_{-h}(s) H_{p,q}(s) \; ds 
			+ \int_{-h}^z \pz \phi_{-h}(z) \phi_{h}(s) H_{p,q}(s) \; ds 
			\\
			= \frac{a_{p,q} c_{p,q}}{2D} \int_{z}^h e^{\bgamma_{p,q}(z+s)} H_{p,q}(s) \; ds 
			+ \frac{a_{p,q} d_{p,q}}{2D} \int_{z}^h e^{-\bgamma_{p,q}(s-z-2h)} H_{p,q}(s) \; ds
			\\
			- \frac{b_{p,q} c_{p,q}}{2D} \int_z^h e^{\bgamma_{p,q}(s-z-2h)} H_{p,q}(s) \; ds 
			- \frac{b_{p,q} d_{p,q}}{2D} \int_z^h e^{-\bgamma_{p,q}(z+s)} H_{p,q}(s) \; ds 
			\\
			+ \frac{a_{p,q} c_{p,q}}{2D} \int_{-h}^z e^{\bgamma_{p,q}(z+s)} H_{p,q}(s) \; ds 
			+ \frac{b_{p,q} c_{p,q}}{2D} \int_{-h}^z e^{-\bgamma_{p,q}(s-z+2h)} H_{p,q}(s) \; ds
			\\
			- \frac{a_{p,q} d_{p,q}}{2D} \int_{-h}^z e^{\bgamma_{p,q}(s-z+2h)} H_{p,q}(s) \; ds 
			- \frac{b_{p,q} d_{p,q}}{2D} \int_{-h}^z e^{-\bgamma_{p,q}(z+s)} H_{p,q}(s) \; ds
			\\
			= \frac{a_{p,q} c_{p,q}}{2D} \int_{-h}^h e^{\bgamma_{p,q}(z+s)} H_{p,q}(s) \; ds 
			- \frac{b_{p,q} d_{p,q}}{2D} \int_{-h}^h e^{-\bgamma_{p,q}(z+s)} H_{p,q}(s) \; ds
			\\
			+ \frac{a_{p,q} d_{p,q}}{2D} \int_{z}^h e^{-\bgamma_{p,q}(s-z-2h)} H_{p,q}(s) \; ds 
			+ \frac{b_{p,q} c_{p,q}}{2D} \int_{-h}^z e^{-\bgamma_{p,q}(s-z+2h)} H_{p,q}(s) \; ds
			\\
			- \frac{b_{p,q} c_{p,q}}{2D} \int_{z}^h e^{\bgamma_{p,q}(s-z-2h)} H_{p,q}(s) \; ds 
			- \frac{a_{p,q} d_{p,q}}{2D} \int_{-h}^z e^{\bgamma_{p,q}(s-z+2h)} H_{p,q}(s) \; ds.
		\end{multline}
		From this we see that the estimates of $T_5, \ldots, T_{10}$ can be used to estimate 
		$I_3$. In fact, by using the triangle inequality, we obtain
		\begin{align*}
			I_3 & := \frac{\Abs{\alpha_p}}{\bgamma_{p,q}^2}
			\left( \int_{-h}^h \Abs{ \pz I_h[H_{p,q}] + \pz I_{-h}[H_{p,q}] }^2 \; dz \right)^{1/2} \\
			& \leq \frac{\Abs{\alpha_p}}{\bgamma_{p,q}}
			\left( T_5 + T_6 + T_7 + T_8 + T_9 + T_{10} \right) 
			\leq C \frac{\Norm{H_{p,q}}{L^2}}{\bgamma_{p,q}^2},
		\end{align*}
		where
		\bes
		\lim \limits_{\Abs{(p,q)} \rightarrow \infty}
		\frac{\Abs{\alpha_p}}{\bgamma_{p,q}} = 1,
		\ees
		was used to obtain the last inequality above. From all of this we find
		\begin{multline}
			\label{Eqn:vz:Est}
			\Norm{\hat{v}^z_{p,q}}{L^2} 
			\leq \frac{1}{\bgamma_{p,q}^2} \Norm{\hat{F}^z_{p,q}}{L^2} 
			+ C \frac{\Abs{\hat{Q}_{p,q}}}{\sqrt{\bgamma_{p,q}} \Abs{d_{p,q}}} 
			+ C \frac{\Abs{\hat{R}_{p,q}}}{\sqrt{\bgamma_{p,q}} \Abs{a_{p,q}}} 
			\\
			+ C \frac{\Norm{H_{p,q}}{L^2}}{\bgamma_{p,q}^2} 
			+ C \frac{\Norm{L_{p,q}}{L^2}}{\bgamma_{p,q}^2}.
		\end{multline}
		From the estimates \eqref{Eqn:vx:Est}, \eqref{Eqn:vy:Est}, and \eqref{Eqn:vz:Est}
		we find
		\begin{align*}
			\Norm{\hat{v}_{p,q}}{L^2}^2 
			& = \Norm{\hat{v}^x_{p,q}}{L^2}^2 
			+ \Norm{\hat{v}^y_{p,q}}{L^2}^2
			+ \Norm{\hat{v}^z_{p,q}}{L^2}^2 \\
			& \leq C \left( \frac{\Abs{\hat{Q}_{p,q}}^2}{\bgamma_{p,q}^3} 
			+ \frac{\Abs{\hat{R}_{p,q}}^2}{\bgamma_{p,q}^3} \right)
			+ C \left( \frac{\Norm{H_{p,q}}{L^2}^2}{\bgamma_{p,q}^4} 
			+ \frac{\Norm{L_{p,q}}{L^2}^2}{\bgamma_{p,q}^4} \right)
			+ \frac{\Norm{\hat{F}^z_{p,q}}{L^2}^2}{\bgamma_{p,q}^4}.
		\end{align*}
		In addition, from \eqref{Eqn:H} and \eqref{Eqn:L}, we have
		\begin{align*}
			H_{p,q}(z) & = \frac{-i \alpha_p}{k_0^2 \bepsilon}
			(i \alpha_p \hat{F}^x_{p,q} + i \beta_q \hat{F}^y_{p,q} 
			+ \pz \hat{F}^z_{p,q}) - \hat{F}^x_{p,q}, \\
			L_{p,q}(z) & = \frac{-i \beta_q}{k_0^2 \bepsilon}
			(i \alpha_p \hat{F}^x_{p,q} + i \beta_q \hat{F}^y_{p,q} 
			+ \pz \hat{F}^z_{p,q}) - \hat{F}^y_{p,q},
		\end{align*}
		so that
		\begin{align*}
			\Norm{\hat{v}_{p,q}}{L^2}^2
			& \leq C \left( \frac{\Abs{\hat{Q}_{p,q}}^2}{\bgamma_{p,q}^3}
			+ \frac{\Abs{\hat{R}_{p,q}}^2}{\bgamma_{p,q}^3} \right)
			+ C \frac{\Abs{\alpha_p}^2}{\bgamma_{p,q}^4}
			\Norm{ i \alpha_p \hat{F}^x_{p,q} + i \beta_q \hat{F}^y_{p,q} 
				+ \pz \hat{F}^x_{p,q} }{L^2}^2 \\
			& \quad
			+ C \frac{\Norm{\hat{F}^x_{p,q}}{L^2}^2
				+ \Norm{\hat{F}^y_{p,q}}{L^2}^2
				+ \Norm{\hat{F}^z_{p,q}}{L^2}^2}
			{\bgamma_{p,q}^4}.
		\end{align*}
		Therefore, we obtain
		\begin{align}
			\sump \sumq \Angle{(p,q)}^4 \Norm{\hat{v}_{p,q}}{L^2}^2
			& \leq C \left( \SobNorm{Q}{1/2}^2 + \SobNorm{R}{1/2}^2 
			+ \sump \sumq \Angle{(p,q)}^0 \Norm{\hat{F}_{p,q}}{L^2}^2 
			\right. \notag \\ 
			& \quad \left.
			+ \sump \sumq \Angle{(p,q)}^{2}
			\Norm{ i \alpha_p \hat{F}^x_{p,q} + i \beta_q \hat{F}^y_{p,q} 
				+ \pz \hat{F}^z_{p,q} }{L^2}^2 \right).
			\label{Eqn:Est:v}
		\end{align}
		
		Next, we estimate
		\bes
		\sump \sumq \Angle{(p,q)}^2 \Norm{\pz \hat{v}_{p,q}}{L^2}^2.
		\ees
		For any integer $j \geq 0$, we have 
		\begin{align*}
			\pz^j \phi_h & = \bgamma_{p,q}^j \left( \frac{a_{p,q}}{D} e^{\bgamma_{p,q}(z+h)} 
			+ (-1)^j \frac{b_{p,q}}{D} e^{-\bgamma_{p,q}(z+h)} \right), \\
			\pz^j \phi_{-h} & = \bgamma_{p,q}^j \left( \frac{c_{p,q}}{2 \bgamma_{p,q}} e^{\bgamma_{p,q}(z-h)}
			+ (-1)^j \frac{d_{p,q}}{2 \bgamma_{p,q}} e^{-\bgamma_{p,q}(z-h)} \right),
		\end{align*}
		and, from the Helmholtz equation,
		\bes
		\pz^2 \phi_{h} = \bgamma_{p,q}^2 \phi_{h},
		\quad
		\pz^2 \phi_{-h} = \bgamma_{p,q}^2 \phi_{-h}.
		\ees
		From Lemma~\ref{Lemma:BVP:vx:vy}, we also notice that, for $j \geq 0$,
		\begin{align*}
			\pz^j \hat{v}^x_{p,q}(z) & = -\hat{Q}^x_{p,q} \pz^j \phi_{h}(z;p,q) 
			- \hat{R}^x_{p,q} \pz^j \phi_{-h}(z;p,q) \frac{2 \bgamma_{p,q}}{D} \\
			& \quad
			- \pz^j (I_h[H_{p,q}](z) + I_{-h}[H_{p,q}](z)).
		\end{align*}
		Therefore, from \eqref{Eqn:pzIh}, we obtain
		\begin{align*}
			\pz \hat{v}^x_{p,q} & = \bgamma_{p,q} \left[ -\hat{Q}^x_{p,q} \frac{a_{p,q}}{D}
			e^{\bgamma_{p,q}(z+h)} + \hat{Q}^x_{p,q} \frac{b_{p,q}}{D} e^{-\bgamma_{p,q}(z+h)} 
			\right. \\ & \quad \left.
			- \hat{R}^x_{p,q} \frac{c_{p,q}}{D} e^{\bgamma_{p,q}(z-h)}
			+ \hat{R}^x_{p,q} \frac{d_{p,q}}{D} e^{-\bgamma_{p,q}(z-h)} \right] \\
			& \quad
			- \frac{a_{p,q} c_{p,q}}{2D} \int_{-h}^h e^{\bgamma_{p,q}(z+s)} H_{p,q}(s) \; ds 
			+ \frac{b_{p,q} d_{p,q}}{2D} \int_{-h}^h e^{-\bgamma_{p,q}(z+s)} H_{p,q}(s) \; ds \\
			& \quad 
			+ \frac{a_{p,q} d_{p,q}}{2D} \int_{z}^h e^{\bgamma_{p,q} (z-s+2h)} H_{p,q}(s) \; ds 
			+ \frac{b_{p,q} c_{p,q}}{2D} \int_{-h}^z e^{\bgamma_{p,q}(z-s-2h)} H_{p,q}(s) \; ds \\
			& \quad 
			- \frac{b_{p,q} c_{p,q}}{2D} \int_{z}^h e^{-\bgamma_{p,q}(z-s+2h)} H_{p,q}(s) \; ds 
			- \frac{a_{p,q} d_{p,q}}{2D} \int_{-h}^z e^{-\bgamma_{p,q}(z-s-2h)} H_{p,q}(s) \; ds.
		\end{align*}
		We observe that the explicit form of $\pz \hat{v}^x_{p,q}$ is quite similar to that of
		$\hat{v}^x_{p,q}$ so we can use the estimates for $T_1, \ldots, T_{10}$ to estimate 
		$\pz \hat{v}^x_{p,q}$. 
		In particular we find
		\begin{align*}
			\Norm{\pz \hat{v}^x_{p,q}}{L^2}
			& \leq \bgamma_{p,q}(T_1 + \ldots + T_{10}) \\
			& \leq C \left( \frac{\Abs{\hat{Q}^x_{p,q}}}{\sqrt{\bgamma_{p,q}}} 
			+ \frac{\Abs{\hat{R}^x_{p,q}}}{\sqrt{\bgamma_{p,q}}} \right) 
			+ C \frac{\Norm{H_{p,q}}{L^2}}{\bgamma_{p,q}}.
		\end{align*}
		Similarly, we also obtain
		\begin{align*}
			\Norm{\pz \hat{v}^y_{p,q}}{L^2}
			& \leq \bgamma_{p,q} (T_1 + \ldots + T_{10}) \\
			& \leq C \left( \frac{\Abs{\hat{Q}^y_{p,q}}}{\sqrt{\bgamma_{p,q}}}
			+ \frac{\Abs{\hat{R}^y_{p,q}}}{\sqrt{\bgamma_{p,q}}} \right)
			+ C \frac{\Norm{L_{p,q}}{L^2}}{\bgamma_{p,q}}.
		\end{align*}
		Next, by using \eqref{Eqn:BVP:vx} and \eqref{Eqn:BVP:vy}, we have
		\be
		\label{Eqn:vxzz:vyzz}
		\pz^2 \hat{v}^x_{p,q} = \bgamma_{p,q}^2 \hat{v}^x_{p,q} 
		+ H_{p,q}(z),
		\quad
		\pz^2 \hat{v}^y_{p,q} = \bgamma_{p,q}^2 \hat{v}^y_{p,q} 
		+ L_{p,q}(z).
		\ee
		Therefore,
		\begin{align*}
			\pz \hat{v}^z_{p,q} 
			& = \frac{1}{\bgamma_{p,q}^2} \pz \hat{F}^z_{p,q} 
			- \frac{i \alpha_p}{\bgamma_{p,q}^2} \pz^2 \hat{v}^x_{p,q} 
			- \frac{i \beta_q}{\bgamma_{p,q}^2} \pz^2 \hat{v}^y_{p,q} \\
			& = \frac{1}{\bgamma_{p,q}^2} \pz \hat{F}^z_{p,q} 
			- i \alpha_p \hat{v}^x_{p,q} 
			- \frac{i\alpha_p}{\bgamma_{p,q}^2} H_{p,q}(z) 
			- i \beta_q \hat{v}^y_{p,q} 
			- \frac{i \beta_q}{\bgamma_{p,q}^2} L_{p,q}(z).
		\end{align*}
		Thus, we obtain
		\begin{align*}
			\Norm{\pz \hat{v}^z_{p,q}}{L^2} 
			& \leq \frac{1}{\bgamma_{p,q}^2} \Norm{\pz \hat{F}^z_{p,q}}{L^2} 
			+ \Abs{\alpha_p} \Norm{\hat{v}^x_{p,q}}{L^2} 
			+ \Abs{\beta_q} \Norm{\hat{v}^y_{p,q}}{L^2} \\
			& \quad
			+ \frac{\Abs{\alpha_p}}{\bgamma_{p,q}^2} \Norm{H_{p,q}}{L^2}
			+ \frac{\Abs{\beta_q}}{\bgamma_{p,q}^2} \Norm{L_{p,q}}{L^2}.
		\end{align*}
		So, recalling the estimates of $\hat{v}^x_{p,q}$ and $\hat{v}^y_{2n}$ in
		\eqref{Eqn:vx:Est} and \eqref{Eqn:vy:Est}, we obtain
		\bes
		\Norm{\pz \hat{v}^z_{p,q}}{L^2} 
		\leq \frac{1}{\bgamma_{p,q}^2} \Norm{\pz \hat{F}^z_{p,q}}{L^2} 
		+ C \left( \frac{\Abs{\hat{Q}_{p,q}}}{\bgamma_{p,q}} 
		+ \frac{\Abs{\hat{R}_{p,q}}}{\sqrt{\bgamma_{p,q}}} \right)
		+ C \left( \frac{\Norm{H_{p,q}}{L^2}}{\bgamma_{p,q}} 
		+ \frac{\Norm{L_{p,q}}{L^2}}{\bgamma_{p,q}} \right).
		\ees
		Combining the estimates of $\pz \hat{v}^x_{p,q}$, $\pz \hat{v}^y_{p,q}$,
		and $\pz \hat{v}^z_{p,q}$ gives 
		\begin{align*}
			\Norm{\pz \hat{v}_{p,q}}{L^2}^2
			& = \Norm{\pz \hat{v}^x_{p,q}}{L^2}^2 
			+ \Norm{\pz \hat{v}^y_{p,q}}{L^2}^2
			+ \Norm{\pz \hat{v}^z_{p,q}}{L^2}^2 \\
			& \leq C \left( \frac{\Abs{\hat{Q}_{p,q}}^2}{\bgamma_{p,q}} 
			+ \frac{\Abs{\hat{R}_{p,q}}^2}{\bgamma_{p,q}} \right) 
			+ C \left( \frac{\Norm{H_{p,q}}{L^2}^2}{\bgamma_{p,q}^2} 
			+ \frac{\Norm{L_{p,q}}{L^2}^2}{\bgamma_{p,q}^2} \right) 
			\\ & \quad
			+ C \frac{1}{\bgamma_{p,q}^4} \Norm{\pz \hat{F}^z_{p,q}}{L^2}^2,
		\end{align*}
		which results in
		\begin{align}
			\sump \sumq \Angle{(p,q)}^2 \Norm{\pz \hat{v}_{p,q}}{L^2}^2
			& \leq C \left( \SobNorm{Q}{1/2}^2 + \SobNorm{R}{1/2}^2 \right) \notag \\
			& \quad
			+ C \sump \sumq ( \Norm{H_{p,q}}{L^2}^2 + \Norm{L_{p,q}}{L^2}^2 ) \notag \\
			& \quad 
			+ C \sump \sumq \Angle{(p,q)}^{-2} \Norm{\pz \hat{F}^z_{p,q}}{L^2}^2 \notag \\
			& \leq C \left( \SobNorm{Q}{1/2}^2 + \SobNorm{R}{1/2}^2 
			+ \sump \sumq \Angle{(p,q)}^{0} \Norm{\hat{F}_{p,q}}{L^2}^2 
			\right. \notag \\
			& \quad
			+ \sump \sumq \Angle{(p,q)}^{2} \Norm{ i \alpha_p \hat{F}^x_{p,q} 
				+ i \beta_q \hat{F}^y_{p,q} + \pz \hat{F}^z_{p,q} }{L^2}^2 \notag \\
			& \quad \left.
			+ \sump \sumq \Angle{(p,q)}^{-2} \Norm{\pz \hat{F}^z_{p,q}}{L^2}^2 \right).
			\label{Eqn:Est:vz}
		\end{align}
		
		We conclude with the estimate of the final sum
		\bes
		\sump \sumq \Norm{\pz^2 \hat{v}_{p,q}}{L^2}^2. 
		\ees
		From \eqref{Eqn:vxzz:vyzz} we have
		\begin{align*}
			\Norm{\pz^2 \hat{v}^x_{p,q}}{L^2} 
			& \leq \bgamma_{p,q}^2 \Norm{\hat{v}^x_{p,q}}{L^2} 
			+ \Norm{H_{p,q}}{L^2}, \\
			\Norm{\pz^2 \hat{v}^y_{p,q}}{L^2} 
			& \leq \bgamma_{p,q}^2 \Norm{\hat{v}^y_{p,q}}{L^2} 
			+ \Norm{L_{p,q}}{L^2},
		\end{align*}
		and the estimates of $\hat{v}^x_{p,q}$ and $\hat{v}^y_{p,q}$
		give the following
		\begin{align*}
			\Norm{\pz^2 \hat{v}^x_{p,q}}{L^2} 
			& \leq C \sqrt{\bgamma_{p,q}} \left( \Abs{\hat{Q}^x_{p,q}}
			+ \Abs{\hat{R}^x_{p,q}} \right) + C \Norm{H_{p,q}}{L^2}, \\
			\Norm{\pz^2 \hat{v}^y_{p,q}}{L^2} 
			& \leq C \sqrt{\bgamma_{p,q}} \left( \Abs{\hat{Q}^y_{p,q}} 
			+ \Abs{\hat{R}^y_{p,q}} \right) + C \Norm{L_{p,q}}{L^2}.
		\end{align*}
		To conclude, we require an estimate of $\pz^2 \hat{v}^z_{p,q}$,
		and, by using \eqref{Eqn:vxzz:vyzz} again, we have
		\begin{align*}
			\pz^2 \hat{v}^z_{p,q} 
			& = \frac{1}{\bgamma_{p,q}^2} \pz^2 \hat{F}^z_{p,q} 
			- \frac{i \alpha_p}{\bgamma_{p,q}^2} \pz^3 \hat{v}^x_{p,q} 
			- \frac{i \beta_q}{\bgamma_{p,q}^2} \pz^3 \hat{v}^y_{p,q} \\
			& =\frac{1}{\bgamma_{p,q}^2} \pz^2 \hat{F}^z_{p,q} 
			- i \alpha_p \pz \hat{v}^x_{p,q} 
			- i \beta_q \pz \hat{v}^y_{p,q} 
			- \frac{i \alpha_p}{\bgamma_{p,q}^2} \pz H_{p,q}(z) 
			- \frac{i \beta_q}{\bgamma_{p,q}^2} \pz L_{p,q}(z).
		\end{align*}
		
		From the estimates of $\pz \hat{v}^x_{p,q}$ and $\pz \hat{v}^y_{p,q}$
		we can bound $\pz^2 \hat{v}^z_{p,q}$. More specifically,
		\begin{align*}
			\Norm{ \pz^2 \hat{v}^z_{p,q} }{L^2} 
			& \leq \frac{1}{\bgamma_{p,q}^2} \Norm{\pz^2 \hat{F}^z_{p,q}}{L^2}
			+ C \left( \sqrt{\bgamma_{p,q}} \Abs{\hat{Q}_{p,q}} 
			+ \sqrt{\bgamma_{p,q}} \Abs{\hat{R}_{p,q}} \right) \\
			& \quad
			+ C \left( \Norm{H_{p,q}}{L^2} + \Norm{L_{p,q}}{L^2} \right)
			+ C \left( \frac{\Norm{\pz H_{p,q}}{L^2}}{\bgamma_{p,q}}
			+ \frac{\Norm{\pz L_{p,q}}{L^2}}{\bgamma_{p,q}} \right).
		\end{align*}
		From this we obtain
		\begin{align*}
			\Norm{\pz^2 \hat{v}_{p,q}}{L^2}^2
			& \leq \frac{C}{\bgamma_{p,q}^4} \Norm{\pz^2 \hat{F}^z_{p,q}}{L^2}^2
			+ C \left( \bgamma_{p,q} \Abs{\hat{Q}_{p,q}}^2 
			+ \bgamma_{p,q} \Abs{\hat{R}_{p,q}}^2 \right) \\
			& \quad
			+ C \left( \Norm{H_{p,q}}{L^2}^2 + \Norm{L_{p,q}}{L^2}^2 \right)
			+ C \left( \frac{\Norm{\pz H_{p,q}}{L^2}^2}{\bgamma_{p,q}^2} 
			+ \frac{\Norm{\pz L_{p,q}}{L^2}^2}{\bgamma_{p,q}^2} \right) \\
			& \leq \frac{C}{\bgamma_{p,q}^4} \Norm{\pz \hat{F}^z_{p,q}}{L^2}^2
			+ C \left( \bgamma_{p,q} \Abs{\hat{Q}_{p,q}}^2 
			+ \bgamma_{p,q} \Abs{\hat{R}_{p,q}}^2 \right) \\
			& \quad
			+ C \Abs{\alpha_p}^2 \Norm{ i\alpha_p \hat{F}^x_{p,q} 
				+ i \beta_q \hat{F}^y_{p,q} + \pz \hat{F}^z_{p,q}}{L^2}^2 \\
			& \quad
			+ C \left( \Norm{\hat{F}^x_{p,q}}{L^2}^2 + \Norm{\hat{F}^y_{p,q}}{L^2}^2 \right) \\
			& \quad
			+  C \Norm{ i \alpha_p \pz \hat{F}^x_{p,q} + i \beta_q \pz \hat{F}^y_{p,q}
				+ \pz^2 \hat{F}^z_{p,q} }{L^2}^2 \\
			& \quad
			+ \frac{C}{\bgamma_{p,q}^2} \left( \Norm{\pz \hat{F}^x_{p,q}}{L^2}^2
			+ \Norm{\pz \hat{F}^y_{p,q}}{L^2}^2 \right),
		\end{align*}
		which produces
		\begin{align}
			\sump \sumq \Angle{(p,q)}^0 \Norm{\pz^2 \hat{v}_{p,q}}{L^2}^2
			& \leq C \left( \SobNorm{Q}{1/2}^2 + \Norm{R}{L^2}^2 
			+ \sump \sumq \Norm{\hat{F}_{p,q}}{L^2}^2 
			\right. \notag \\
			& \quad
			+ \sump \sumq \Angle{(p,q)}^2 \Norm{ \hat{F}^x_{p,q} }{L^2}^2 
			+ \Norm{ \hat{F}^y_{p,q} }{L^2}^2 \notag \\
			& \quad
			+ \sump \sumq \Norm{ i \alpha_p \pz \hat{F}^x_{p,q} 
				+ i \beta_q \pz \hat{F}^y_{p,q} 
				+ \pz^2 \hat{F}^z_{p,q} }{L^2}^2 \notag \\
			& \quad
			+ \sump \sumq \Angle{(p,q)}^{-2} \Norm{ \pz \hat{F}_{p,q} }{L^2}^2 \notag \\
			& \quad \left.
			+ \sump \sumq \Angle{(p,q)}^{-4} \Norm{ \pz^2 \hat{F}^z_{p,q} }{L^2}^2 \right).
			\label{Eqn:Est:vzz}
		\end{align}
		Moreover,
		\bes
		\Div{F} = \sump \sumq (i \alpha_p \hat{F}^x_{p,q} + i \beta_q \hat{F}^y_{p,q}
		+ \pz \hat{F}^z_{p,q}) e^{i \alpha_p x + i \beta_q y},
		\ees
		so
		\begin{align*}
			\SobNorm{\Div{F}}{1}^2 
			& = \sump \sumq \Angle{(p,q)}^2 \Norm{ i \alpha_p \hat{F}^x_{p,q} 
				+ i \beta_q \hat{F}^y_{p,q} + \pz \hat{F}^z_{p,q} }{L^2}^2 \\
			& \quad
			+ \sump \sumq \Angle{(p,q)}^0 \Norm{ \pz ( i \alpha_p \hat{F}^x_{p,q} 
				+ i \beta_q \hat{F}^y_{p,q} + \pz \hat{F}^z_{p,q} ) }{L^2}^2.
		\end{align*}
		Finally, we combine all the estimates from \eqref{Eqn:Est:v}, \eqref{Eqn:Est:vz},
		and \eqref{Eqn:Est:vzz} to obtain
		\begin{align*}
			\SobNorm{v}{2}^2 
			& \leq C \left( \SobNorm{Q}{1/2}^2 + \SobNorm{R}{1/2}^2 \right) \\
			& \quad
			+ C \left( \sump \sumq \Angle{(p,q)}^2 \Norm{ i \alpha_p \hat{F}^x_{p,q} 
				+ i \beta_q \hat{F}^y_{p,q} + \pz \hat{F}^z_{p,q} }{L^2}^2 
			\right. \\
			& \quad \left.
			+ \sump \sumq \Norm{ \pz ( i \alpha_p \hat{F}^x_{p,q} 
				+ i \beta_q \hat{F}^y_{p,q} + \pz \hat{F}^z_{p,q} ) }{L^2}^2 \right) \\
			& \quad 
			+ C \sump \sumq \Angle{(p,q)}^{-4} \left( \Angle{(p,q)}^4 \Norm{\hat{F}_{p,q}}{L^2}^2
			+ \Angle{(p,q)}^2 \Norm{\pz \hat{F}_{p,q}}{L^2}^2 
			+ \Norm{\pz^2 \hat{F}_{p,q} }{L^2}^2 \right) \\
			& \leq 
			C \left( \SobNorm{Q}{1/2}^2 + \SobNorm{R}{1/2}^2 
			+ \SobNorm{\Div{F}}{1}^2 + \SobNorm{F}{2}^2 \right).
		\end{align*}
		With this the proof is complete.
	\end{proof}
	
	%
	%
	
	\section{Proof of an Elliptic Estimate: Joint Analyticity}
	\label{Sec:JointAnal:Proof}
	
	In this appendix we establish the elliptic estimate,
	Theorem~\ref{Thm:JointAnal:EllEst},
	required to prove joint analyticity of the scattered field using the
	inductive proof described in Section~\ref{Sec:JointAnal}. As our
	proof is inductive in the vertical, $z$, derivative, we begin with the
	following lateral regularity result.
	\begin{theorem}
		\label{Thm:JointAnal:EllEst:Zero}
		Given any integer $m \geq 0$, if
		$(\omega,\bepsilon) \in \cP$,
		$F(x,y,z) \in C^{\omega}_m(S_v)$
		such that
		\bes
		\max \left\{
		\SobNorm{ \frac{\px^r \py^t}{(r+t)!} F }{0},
		\SobNorm{ \frac{\px^r \py^t}{(r+t)!} \Div{F} }{1}
		\right\}
		\leq C_F \frac{\eta^r}{(r+1)^2} \frac{\theta^t}{(t+1)^2},
		\ees
		for all $r, t \geq 0$ and for some constants 
		$C_F, \eta, \theta > 0$, and
		$Q, R \in C^{\omega}_m(\Gamma))$ satisfying 
		\begin{align*}
			& \SobNorm{ \frac{\px^r \py^t}{(r+t)!} Q }{1/2} 
			\leq C_Q \frac{\eta^r}{(r+1)^2} \frac{\theta^t}{(t+1)^2}, \\
			& \SobNorm{\frac{\px^r \py^t}{(r+t)!} R }{1/2} 
			\leq C_R \frac{\eta^r}{(r+1)^2} \frac{\theta^t}{(t+1)^2},
		\end{align*}
		for all $r, t \geq 0$ and some constants $C_R, C_Q > 0$.
		Then, there exists a unique solution
		$v \in C^{\omega}_m(S_v)$ of 
		\bse
		\label{Eqn:JointAnal:EllEst}
		\begin{align}
			& \cL_0 v = F, && \text{in $S_v$}, \\
			& -\pz v - T_u [v] = Q, && \text{at $\Gamma_h$}, \\
			& \pz v - T_w[v] = R, && \text{at $\Gamma_{-h}$}, \\
			& v(x+d_x,y+d_y, z) = e^{i \alpha d_x + i \beta d_y} v(x,y,z),
		\end{align}
		\ese
		satisfying 
		\bes
		\SobNorm{ \frac{\px^r \py^t}{(r+t)!} v }{2} 
		\leq \underline{C}_e \frac{\eta^r}{(r+1)^2} \frac{\theta^t}{(t+1)^2},
		\ees
		for all $r, t \geq 0$ where
		\bes
		\underline{C}_e = C_e (C_F + C_Q + C_R) > 0,
		\ees
		and $\overline{C}(h) > 0$ is a universal constant.
	\end{theorem}
	\begin{proof}
		Applying the operator $\p_x^r \p_y^t/(r+t)!$ to 
		\eqref{Eqn:JointAnal:EllEst} we obtain
		\begin{align*}
			& \cL_0 \left[ \frac{\px^r \py^t}{(r+t)!} v \right] 
			= \frac{\px^r \py^t}{(r+t)!} F, && \text{in $S_v$}, \\
			& -\pz \left[ \frac{\px^r \py^t}{(r+t)!} v \right]
			- T_u \left[ \frac{\px^r \py^t}{(r+t)!} v \right] 
			= \frac{\px^r \py^t}{(r+t)!} Q, && \text{at $\Gamma_h$}, \\
			& \pz \left[ \frac{\px^r \py^t}{(r+t)!} v \right]
			- T_w \left[ \frac{\px^r \py^t}{(r+t)!} v \right] 
			= \frac{\px^r \py^t}{(r+t)!} R, && \text{at $\Gamma_{-h}$}, \\
			& \frac{\px^r \py^t}{(r+t)!} v(x+d_x,y+d_y, z) 
			= e^{i \alpha d_x + i \beta d_y} \frac{\px^r \py^t}{(r+t)!} v(x,y,z).
		\end{align*}
		Applying Theorem~\ref{Thm:EllEst} to this system and using the
		hypothesis on $F$, $Q$, and $R$, we obtain
		\begin{align*}
			\SobNorm{ \frac{\px^r \py^t}{(r+t)!} v }{2} 
			& \leq C_e \left( \SobNorm{ \frac{\px^r \py^t}{(r+t)!} F }{0}
			+ \SobNorm{ \Div{ \frac{\px^r \py^t}{(r+t)!} F } }{1}
			\right. \\ & \quad \left.
			+ \SobNorm{ \frac{\px^r \py^t}{(r+t)!} Q }{1/2}
			+ \SobNorm{ \frac{\px^r \py^t}{(r+t)!} R }{1/2} \right) \\
			& \leq C_e (C_F + C_Q + C_R) 
			\frac{\eta^r}{(r+1)^2} \frac{\theta^t}{(t+1)^2},
		\end{align*}
		and we are done by choosing  $\underline{C_e} = C_e (C_F+C_Q+C_R)$.
	\end{proof}
	
	We are now in a position to establish our main result.			
	\begin{proof}{[Theorem~\ref{Thm:JointAnal:EllEst}]}.
		We prove \eqref{Eqn:JointAnal:EllEst:Estimate} by induction in $s$,
		and the case $s = 0$ is verified by the previous result,
		Theorem~\ref{Thm:JointAnal:EllEst:Zero}.
		We now assume that 
		\bes
		\SobNorm{ \frac{\px^r \py^t \pz^s}{(r+t+s)!} v }{2}
		\leq \underline{C_e} \frac{\eta^r}{(r+1)^2} 
		\frac{\theta^t}{(t+1)^2} \frac{\zeta^s}{(s+1)^2}
		\quad
		\forall\ s < \overline{s},
		\quad
		\forall\ r, t \geq 0,
		\ees
		which implies that, for $m \in \{ x, y, z \}$,
		\bes
		\SobNorm{ \frac{\px^r \py^t \pz^s}{(r+t+s)!} v^m }{2}
		\leq \underline{C_e} \frac{\eta^r}{(r+1)^2} 
		\frac{\theta^t}{(t+1)^2} \frac{\zeta^s}{(s+1)^2}
		\quad
		\forall\ s < \overline{s},
		\quad
		\forall\ r, t \geq 0.
		\ees
		We now examine
		\begin{align*}
			\SobNorm{ \frac{\px^r \py^t \pz^{\overline{s}}}
				{(r+t+\overline{s})!} v }{2}
			& \leq \SobNorm{ \frac{\px^r \py^t \pz^{\overline{s}}}
				{(r+t+\overline{s})!} v }{1}
			+ \SobNorm{ \frac{\px^r \py^t \pz^{\overline{s}}}
				{(r+t+\overline{s})!} \px v }{1} \\
			& \quad
			+ \SobNorm{ \frac{\px^r \py^t \pz^{\overline{s}}}
				{(r+t+\overline{s})!} \py v }{1}
			+ \SobNorm{ \frac{\px^r \py^t \pz^{\overline{s}}}
				{(r+t+\overline{s})!} \pz v }{1} \\
			& \leq \SobNorm{ \frac{\px^r \py^t \pz^{\overline{s}-1}}
				{(r+t+\overline{s})!} v }{2}
			+ \SobNorm{ \frac{\px^r \py^t \pz^{\overline{s}-1}}
				{(r+t+\overline{s})!} \px v }{2} \\
			& \quad
			+ \SobNorm{ \frac{\px^r \py^t \pz^{\overline{s}-1}}
				{(r+t+\overline{s})!} \py v }{2}
			+ \SobNorm{ \frac{\px^r \py^t \pz^{\overline{s}-1}}
				{(r+t+\overline{s})!} \pz^2 v }{1}.
		\end{align*}
		The first three terms can be bounded by our inductive hypothesis as
		they involve $z$ derivatives of order $\overline{s}-1$. The fourth 
		term we denote $W$ and analyze as follows. Writing
		\bes
		W^m := \SobNorm{ \frac{\px^r \py^t \pz^{\overline{s}-1}}
			{(r+t+\overline{s})!} \pz^2 v^m }{1},
		\quad
		m \in \{ x, y, z \},
		\ees
		it is clear that $W \leq W^x + W^y + W^z$. To estimate $W^x$ and 
		$W^y$ we recall that the inhomogeneous Maxwell equations
		\bes
		-\Laplacian{v} + \Grad{\Div{v}} - k_0^2 \bepsilon v = F,
		\ees
		give
		\bes
		\pz^2 v = -F - k_0^2 \bepsilon v - \px^2 v - \py^2 v
		+ \begin{pmatrix} \px (\px v^x + \py v^y + \pz v^z) \\
			\py (\px v^x + \py v^y + \pz v^z) \\
			\pz (\px v^x + \py v^y + \pz v^z) \end{pmatrix}.
		\ees
		Therefore, we can write $\pz^2 v^x$ and $\pz^2 v^y$ as
		\begin{align*}
			\pz^2 v^x & = -F^x - k_0^2 \bepsilon v^x - \py^2 v^x 
			+ \px \py v^y + \px \pz v^z, \\
			\pz^2 v^y & = -F^y - k_0^2 \bepsilon v^y - \px^2 v^y
			+ \px \py v^x + \py \pz v^z.
		\end{align*}
		With these and the inductive hypothesis on $v$ we obtain
		\begin{align*}
			W^x & \leq \SobNorm{ \frac{\px^r \py^t \pz^{\overline{s}-1}}
				{(r+t+\overline{s})!} F^x }{1} 
			+ k_0^2 \bepsilon \SobNorm{ \frac{\px^r \py^t \pz^{\overline{s}-1}}
				{(r+t+\overline{s})!} v^x }{1} 
			+ \SobNorm{ \frac{\px^r \py^{t+1} \pz^{\overline{s}-1}}
				{(r+t+\overline{s})!} \py v^x }{1} \\
			& \quad
			+ \SobNorm{ \frac{\px^r \py^{t+1} \pz^{\overline{s}-1}}
				{(r+t+\overline{s})!} \px v^y }{1}
			+ \SobNorm{ \frac{\px^{r+1} \py^t \pz^{\overline{s}-1}}
				{(r+t+\overline{s})!} \pz v^z }{1} \\
			& \leq \SobNorm{ \frac{\px^r \py^t \pz^{\overline{s}-1}}
				{(r+t+\overline{s})!} F^x }{2} 
			+ k_0^2 \bepsilon \SobNorm{ \frac{\px^r \py^t \pz^{\overline{s}-1}}
				{(r+t+\overline{s})!} v^x }{2} 
			+ \SobNorm{ \frac{\px^r \py^{t+1} \pz^{\overline{s}-1}}
				{(r+t+\overline{s})!} v^x }{2} \\
			& \quad
			+ \SobNorm{ \frac{\px^r \py^{t+1} \pz^{\overline{s}-1}}
				{(r+t+\overline{s})!} v^y }{2}
			+ \SobNorm{ \frac{\px^{r+1} \py^t \pz^{\overline{s}-1}}
				{(r+t+\overline{s})!} v^z }{2} \\
			& \leq \underline{C_e} \frac{\eta^r}{(r+1)^2} \frac{\theta^t}{(t+1)^2} 
			\frac{\zeta^{\overline{s}}}{(\overline{s}+1)^2} 
			\left( \frac{2 \theta}{\zeta} \frac{(t+1)^2}{(t+2)^2}
			\frac{(\overline{s}+1)^2}{\overline{s}^2} 
			+ \frac{\eta}{\zeta} \frac{(r+1)^2}{(r+2)^2}
			\frac{(\overline{s}+1)^2}{\overline{s}^2} \right) \\
			& \leq \underline{C_e} \frac{\eta^r}{(r+1)^2} \frac{\theta^t}{(t+1)^2}
			\frac{\zeta^{\overline{s}}}{(\overline{s}+1)^2}
			\left( \frac{8 \theta + 4 \eta}{\zeta} \right),
		\end{align*}
		where the last inequality is obtained by using the inequalities
		\bes
		\frac{(\overline{s}+1)^2}{\overline{s}^2} \leq 4,
		\quad
		\frac{(t+1)^2}{(t+2)^2} \leq 1.
		\ees
		In analogous fashion we obtain
		\bes
		W^y \leq \underline{C_e} \frac{\eta^r}{(r+1)^2} \frac{\theta^t}{(t+1)^2}
		\frac{\zeta^{\overline{s}}}{(\overline{s}+1)^2}
		\left( \frac{8 \theta + 4 \eta}{\zeta} \right).
		\ees
		All that remains is to estimate $W^z$ which we accomplish by
		applying the divergence operator to the Maxwell equations
		(noting that $\Div{\Curl{\psi}} = 0$ for any $\psi$) giving
		\bes
		-k_0^2 \bepsilon \Div{v} = \div{F},
		\ees
		which implies that 
		\bes
		\pz v^z = -\frac{1}{k_0^2 \bepsilon} \Div{F} - \px v^x - \py v^y,
		\ees
		and
		\bes
		\pz^2 v^z = -\frac{1}{k_0^2 \bepsilon} \pz \Div{F} - \px \pz v^x - \py \pz v^y.
		\ees
		Therefore, 
		\begin{align*}
			W^z & \leq \frac{1}{k_0^2 \bepsilon}
			\SobNorm{ \frac{\px^r \py^t \pz^{\overline{s}-1}}
				{(r+t+\overline{s})!} \pz \Div{F} }{1}
			+ \SobNorm{ \frac{\px^{r+1} \py^t \pz^{\overline{s}-1}}
				{(r+t+\overline{s})!} \pz v^x }{1} 
			+ \SobNorm{ \frac{\px^r \py^{t+1} \pz^{\overline{s}-1}}
				{(r+t+\overline{s})!} \pz v^y }{1} \\
			& \leq \frac{1}{k_0^2 \bepsilon} 
			\SobNorm{ \frac{\px^r \py^t \pz^{\overline{s}-1}}
				{(r+t+\overline{s})!} \pz \Div{F} }{1} 
			+ \SobNorm{ \frac{\px^{r+1} \py^t \pz^{\overline{s}-1}}
				{(r+t+\overline{s})!} v^x }{2} 
			+ \SobNorm{ \frac{\px^r \py^{t+1} \pz^{\overline{s}-1}}
				{(r+t+\overline{s})!} v^y }{2}.
		\end{align*}
		The first term is readily estimated from the hypothesis of $F$.
		The second and third terms can be addressed by our inductive 
		hypothesis for $v^x$ and $v^y$ as they just involve $z$ 
		derivatives of order $\overline{s}-1$. Thus, we obtain
		\bes
		W^z \leq \underline{C_e} \frac{\eta^r}{(r+1)^2}
		\frac{\theta^t}{(t+1)^2} 
		\frac{\zeta^{\overline{s}}}{(\overline{s}+1)^2}
		\left( \frac{4 \theta + 4 \eta}{\zeta} \right).
		\ees
		Finally, collecting all the estimates for $W^x$, $W^y$, and 
		$W^z$, we arrive at 
		\bes
		W^x + W^y + W^z \leq \underline{C_e} 
		\frac{\eta^r}{(r+1)^2} \frac{\theta^t}{(t+1)^2}
		\frac{\zeta^{\overline{s}}}{(\overline{s}+1)^2}
		\left( \frac{20 \theta + 12 \eta}{\zeta} \right).
		\ees
		The proof is complete by choosing $\zeta \geq 20 \theta + 12 \eta$.
	\end{proof}
	
	%
	%
	
	\bibliographystyle{abbrv}
	\bibliography{nicholls}
	
\end{document}